\mag=\magstep1
\hsize=16.5truecm
\vsize=23.3truecm
\documentstyle{amsppt}
\document
\def\Aut{\operatorname{Aut}}
\def\nsubset{\nsubseteq}
\def\Gal{\operatorname{Gal}}
\def\Coker{\operatorname{Coker}}

\def\Im{\operatorname{Im}}
\def\Cone{\operatorname{Cone}}
\def\Ind{\operatorname{Ind}}
\def\lcm{\operatorname{lcm}}
\topmatter
\title
Convolution theorem for non-degenerate maps and composite singularities
\endtitle
\author
Tomohide Terasoma
\endauthor
\affil
Department of Mathematical Science,\\
University of Tokyo, \\
Komaba, Meguro, Tokyo, 153, Japan \\
Max-Planck-Institut f\"ur Mathematik,\\
Gottfried Claren Strasse 26,\\
D-53225, Bonn, Germany.
\endaffil
\rightheadtext{Convolution Theorem and composite singularities}
\endtopmatter
\heading
\S 0 Introduction
\endheading

Let $n$ be a natural number greater than 1 and $f=f(x) = f(x_1, \dots , x_n)$
be a germ of holomorphic function on $\bold C^n$ at $0$.  Let $\epsilon$
be a sufficiently small number such that $f$ is defined on 
$B(x)= \{ (x_1, \dots ,x_n)\in \bold C^n \mid \mid x_i \mid <\epsilon \}$.
If $\delta'$ is sufficiently small, the restriction of $f:B(x) \to \bold C$ to
$f^{-1}(B^0(s))$, where $B^0(s) = \{ s \in \bold C \mid 
0< \mid s \mid <\delta' \}$ is a fiber bundle over $B^0(s)$.
Therefore for $0<\delta< \delta'$, $\pi_1(B^0(s))=\pi_1(B^0(s), \delta)$ 
acts on $H^i(f^{-1}(\delta ), \bold Q)$.  
From now on, we use the notation $f: B(x) \to B(t)$ for short in this
situation, and the morphism $f$ is called a local morphism.
If $f$ has only isolated singularities, then
$H^i(f^{-1}(\delta ), \bold Q)=0$ if $i \neq 0, n-1$.
In the same way if $m$ is 
a natural number greater than 1 and $g$ is a germ of holomorphic
function on 
$\bold C^m=\{ y=(y_1, \dots , y_m)\}$ at 0, a restriction of $g$ defines a
fiber bundle $g^{-1}(B^0(t)) \to B^0(t)$ over $B^0(t)$.
If $f$ and $g$ has only isolated singularities, $u=f+g$ also has only
isolated singularity.
Let $\gamma_f$ , $\gamma_g$ and $\gamma_{f+g}$ be the action of
the canonical generator of $\pi_1(B^0(s))$, $\pi_1(B^0(t))$
and $\pi_1(B^0(u))$ respectively.
Then the relation between $H^{n+m-1}((f+g)^{-1}(\delta ), \bold Q)$,
$H^{n-1}(f^{-1}(\delta ), \bold Q)$ and 
$H^{m-1}(g^{-1}(\delta ), \bold Q)$ is given by Thom-Sebastiani
theorem [S-T]: There exists an isomorphism
$$
H^{n+m-1}((f+g)^{-1}(\delta ), \bold Q) \simeq
H^{n-1}(f^{-1}(\delta ), \bold Q)\otimes
H^{m-1}(g^{-1}(\delta ), \bold Q).
$$
compatible with the actions
of the monodromies, i.e. $\gamma_{f+g} = \gamma_f \otimes \gamma_g$.

The Hodge theoretic refinement of Thom-Sebastiani theorem
is conjectured by Steenbrink [St], and proved by Varchenko [V],
M.Saito [S] and Denef-Loeser [D-L] in different contexts.
Let $f$ be a germ of holomorphic function $f:B(x) \to B(t)$
such that the restriction of $f$ to $f^{-1}(B^0(t))$ is smooth.
After Steenbrink, 
the space $H^i(f^{-1}(\delta ), \bold Q)$ is equipped with the mixed Hodge
structure for a sufficiently small $\delta$.    
The action of $\gamma_f$ is quasi-unipotent and preserves 
the mixed Hodge structure. 
The minimal natural number $m$ for which the action of $\gamma^m$
is unipotent for $H^i(f^{-1}(\delta ), \bold Q)$ $(0 \leq i \leq 2(n-1))$
is called the exponent
of $f$.  If $m$ is the exponent of $f$, then the action of 
$N = \frac{1}{m}\log (\gamma^m )$ satisfies $NW_k \subset W_{k-2}$.
Therefore the action of $\gamma$ on 
$Gr_{k}^W(H^i(f^{-1}(\delta), \bold Q))$
preserves the Hodge structure and has a finite image.
Therefore $[H^i(f^{-1}(\delta), \bold Q)]
=\sum_k[Gr_k^W(H^i(f^{-1}(\delta), \bold Q))]$
defines an element of the Grothendieck group $K_{MH}(\bold C, \mu_m)$ of
mixed Hodge structures with an action of $\mu_m$.
For an element $V$ in $K_{MH}(\bold C, \mu_m)$, the invariant part
and the complement of invariant part under the action of $\mu_m$
is denoted by $V_1$ and $V_{\neq 1}$, respectively.
The action of $\mu_m$ preserves the space $H^{p,q}$, where
$$
H^{p,q} =
F^pGr_{p+q}^W(H^i(f^{-1}(\delta), \bold C))\cap
\bar F^qGr_{p+q}^W(H^i(f^{-1}(\delta), \bold C))
$$
and $\chi$-part of $H^{p,q}$ is denoted by $H^{p,q}(\chi )$,
for a finite character $\chi$ of $\pi_1(B^0(t), \delta )$.
The set of natural numbers $\{ h^{p,q}(\chi)\}$, where
$h^{p,q}(\chi) = \dim H^{p,q}(\chi)$ is called the spectral pair of $f$
and it is an important invariant for singularities.
The Hodge analog of Thom-Sebastiani theorem is stated as follows.
Let $d_1$ and $d_2$ be the exponent of $f$ and $g$ respectively and
$m$ be the l.c.m. of $d_1$ and $d_2$.
On the mixed Hodge structure  
$V(d_1, d_2) =H^1(\{\sigma^{d_1} + \tau^{d_2} = \delta\}, \bold Q)$,
the group $\mu_{d_1} \times \mu_{d_2} \times \pi_1 (B(u)-\{0\},\delta)$
acts and the exponent is equal to $m = \lcm (d_1, d_2)$.
Therefore $V(d_1, d_2)$ defines an element $[V(d_1, d_2)]$ in 
$K_{MH}(\bold C, \mu_{d_1}\times \mu_{d_2} \times \mu_m)$.
\proclaim{Theorem 0.1 ([V], [S], [D-L])}
Under the above notation, the following equality holds in
the Grothendieck group $K_{MH}(\bold C, \mu_m)$.
$$
\align
&[H^{n+m-1}((f+g)^{-1}(\delta ), \bold Q)] \\
=
& -([V(d_1, d_2)] 
\otimes [H^{n-1}(f^{-1}(\delta), \bold Q)]
\otimes [H^{m-1}(g^{-1}(\delta), \bold Q)])^{\mu_{d_1}\times \mu_{d_2}} \\
 & +[H^{n-1}(f^{-1}(\delta), \bold Q)(\chi_1)]_{\neq 1}
\otimes [H^{m-1}(g^{-1}(\delta), \bold Q)(\chi_2)]_{\neq 1} \\
  & -[H^{n-1}(f^{-1}(\delta), \bold Q)(\chi_1)]
\otimes [H^{m-1}(g^{-1}(\delta), \bold Q)(\chi_2)]. \\
\endalign
$$
\endproclaim
This theorem is generalized by Nemethi-Steenbrink [N-S] as follows.
Let $h(s,t)$ be a function of two variables with an isolated singularity.
The composite $h\circ (f, g)$
of $(f,g):B(x) \times B(y) \to B(s)\times B(t)$
and $h: B(s) \times B(t) \to B(u)$ has singularities and
it is called a composite singularity of two variables.
The spectral pair of the composite singularity of two variables
is computed by those of
$f$, $g$ and $h$ in [N-S].
Note that if $h(s,t) = s+t$, then their theorem is nothing but
Theorem 0.1.
In this paper, we generalize this theorem
for composite singularities of several variables.
We explain the situation of our result.

Let $n$ be a natural number greater than 1 and
$g_i(y_i) =g_i(y_{i1}, \dots, y_{im_i})$ be germs of holomorphic
functions of $m_i$ variables for $1 \leq i \leq n$.
Then $x_i = g_i(y_i)$ defines a local
morphism $g_i :B(y_i)\to B(x_i)$.
Suppose that the restriction of
$g_i$ to $g_i^{-1}(B(x_i)-\{0\})$ is smooth.
Let $f$ be a germ of holomorphic function on $B(x) =\{(x_1, \dots , x_n)
\mid \mid x_i \mid < \delta \ ( 1\leq i \leq n)\}$ 
non-degenerate with respect to
the Newton boundary.
(See Section 1 for the definition of the non-degeneracy condition.)
Roughly speaking, our main theorem (Theorem 3.6.1) gives
the relation between the spectral pair of the composite singularity
$f(g_1(y_1), \dots ,g_n(y_n))$ and those of
$g_1, \dots, g_n$ and $f$. This theorem is called 
the convolution theorem for composite singularities.
Section 3 is devoted to prove this main theorem.

In the paper of [D-L] they proved the motivic version of Thom-Sebastiani
theorem under the assumption of tameness and 
the existence of the resolution of
singularities.  Their theorem, Motivic Thom-Sebastiani theorem implies the 
$l$-adic analog of Thom-Sebastiani theorem.
We prove the $l$-adic analog of the convolution theorem for composite 
singularities
under the assumption of tameness condition in Section 2.
Although the $l$-adic analog is not a deep theorem,
it is worth giving a proof in this paper, 
because from this result, one can
easily expect the formula for the spectral pairs for composite
singularities.

\heading
\S 1 Non degenerate maps and reducing multiplicities
\endheading
\heading
\S 1.1 Toric geometry
\endheading

  Let $k$ be a finite field $\bold F_q$ of $q$ elements or the complex
number field $\bold C$ and $n \geq 1$.  Let $f = f(x_1, \dots , x_n)
\in k [[ x_1, \dots , x_n]]$ be a formal power series on $x_1, \dots , x_n$
over $k$ such that $f(0) = 0$.  We assume that for all $i = 1, \dots, n$,
the coefficient of $x_i^{m_i}$ is not zero
for some $m_i$.  The Newton polygon $\Delta = \Delta (f)$
of $f = \sum_{w \in\bold N^n}a_wx^w$ is defined by the convex 
hull of $\{ w \in \bold N^n \mid a_w \neq 0\}$ and $\bold R_{\geq N}^{(i)}
= \{ (0, \dots ,\overset{i}\to{r}, \dots , 0)\mid r \geq N\}$ for
a sufficiently large number $N$.  First we recall the notion of non-degeneracy
with respect to the Newton boundary.  Let $\sigma$ be a face of
$\Delta$, $L_{\sigma}$ be the affine linear hull of $\sigma$
 in $\bold R^n$ and $L_{\sigma,0}$ be the linear subspace of 
$\bold R^n$ (containing 0) parallel to $L_{\sigma}$.  We define $f_\sigma$
by $f_\sigma = \sum_{w \in \sigma} a_wx^w$.  
If we choose $w_0 \in \sigma$, then $f_{\sigma, 0}= f_{\sigma}/x^{w_0} 
\in k[L_{\sigma, 0} \cap \bold Z^n]$.
If $\sigma$ is a non-compact face, there exists a non-empty subset $J$ of
$[1, n]$ such that
$\sigma$ is equal to $\sum_{j \in J}\bold R_+\bold e_j$ outside of a
compact subset.  The scheme 
$Spec(k[[x_j]]_{j \in J}[x_j^{-1}]_{j\in J})$ is denoted by
$B^0_{\sigma}$.
The function $f$ is said to be non-degenerate with respect to a compact 
face $\sigma$ (resp. a non-compact face $\sigma$)
if $\{ f_{\sigma,0} =0\}$ is a smooth variety on 
$Z_{\sigma} = Spec (k[L_{\sigma, 0} \cap \bold Z^n])$
(resp. $B^0_{\sigma}$).
The series $f$ is said to
be non-degenerate with respect to the Newton boundary if it is non-degenerate
with respect to all the faces $\sigma$ of $\Delta$.
The closed subschemes $\{ x_i =0\}$ and $\{ f=0\}$
in $B(x) = Spec(k[[x_1, \dots ,x_n]])$
are denoted by $Z(x_i)$ and $Z(f)$, respectively.  
If $f$ is non-degenerate with respect to the
Newton boundary, we can construct a 
partial resolution of the singularity
of the hypersurface $Z(f) \cup \cup_{i=1}^n Z(x_i)$, 
using toric geometry.
Let us define the dual fan $F(\Delta )$ of $\Delta$
as the set $\{ \sigma^*(x) \mid x \in \Delta \}$, where
$\sigma^*(x) = 
\{l \in (\bold R^n)^* \mid 
l(\sum_{y \in \Delta}\bold R_+(y-x))\geq 0 \}$.  
As in [Od], we can associate a toric variety $X_F$ for a fan $F$.
Note that
the toric variety $X_{F(\Delta)}$ associated to $F(\Delta)$ is isomorphic
to $\bold P_{\Delta}$ defined in [Dan].  
The dual coordinate of $(\bold R^n)^*$ is written by $w^*=(w_1^*, \dots,
w_n^*)$.
Let $F$ be a simplicial refinement
of $F(\Delta )$ 
such that 
$$
\align
&\text{the restriction of $F$ to the coordinate 
hyperplane $\{ w_i^* =0\}$
is equal to}
\tag{1.1.1} \\
&\text{ that of $F(\Delta)$}. \\
\endalign
$$
The natural morphism from $X_F$ to the affine space 
$\bold A^n(x) = Spec(k[x_1, \dots ,x_n])$ is denoted by 
$\pi : X_F \to \bold A^n(x)$.  The base change of $\pi$ by the morphism
$B(x) = Spec(k[[x_1, \dots , x_n]]) \to \bold A^n(x)$
is denoted by $b(x) : \hat B (x) \to B(x)$.

Now we recall the definitions of quasi-smooth varieties
and quasi-normal crossing divisors. (See [St].)  An integral variety $V$
is called quasi-smooth at a point $p \in V$ if there
exists a local parameter $x_1, \dots ,x_n$ at $p$ with the inclusion
$\Cal O_{V,p} \subset \Cal O_{V,p}[\xi_1, \dots ,\xi_n]$, where
$\xi_1^{d_i} = x_1$,$\dots$ ,$\xi_1^{d_i} = x_1$ such that
(1) $\Cal O_{V, p}[\xi_1, \dots , \xi_n]$ is smooth with
regular parameters $\xi_1, \dots, \xi_n$, and (2)
$\Cal O_{V, p}$ is identified with the invariant ring of
$\Cal O_{V, p}[\xi_1, \dots ,\xi_n]$ under the action of
a subgroup $G$ of $\mu_{d_1} \times\cdots \times \mu_{d_n}$, where the action
of $\mu_{d_1}\times \cdots \times \mu_{d_n}$ is given by 
the multiplication of coordinate $(\xi_1, \dots , \xi_n)$.
If $V$ is quasi-smooth at all $p\in V$, $V$ is called a quasi-smooth variety.
A divisor $D$ in $V$ is called a quasi-normal crossing divisor at $p \in D$,
if it is the image of $\xi_1\cdots \xi_s$ ($0<s\leq n$) under
the projection $Spec(\Cal O_{V,p}[\xi_1, \dots ,\xi_n])
\to Spec(\Cal O_{V,p})$.  The divisor is called quasi-normal corssing,
if it is quasi-normal crossing at all $p \in D$.

Since $F$ is a simplicial fan,  $X_F$ and $\hat B (x)$ are quasi-smooth, and
the reduced part of $(b(x))^{-1}(Z(f)\cup \cup_{i=1}^nZ(x_i))$
is a quasi-normal crossing divisor.

 Now we consider a covering $B(\xi)$ of $B(x)$.  Let $d_1, \dots ,d_n$
be natural numbers such that $\mu_{d_i} \subset k^\times$.
We introduce a system of coordinates $\xi_1, \dots ,\xi_n$ with
$\xi_1^{d_1} = x_1, \dots,\xi_n^{d_n} = x_n$.  Then
$B(\xi) = Spec(k[[\xi_1, \dots , \xi_n]])$ is a finite Galois
covering of $B(x)$ with a Galois group 
$G = \mu_{d_1} \times \cdots \times \mu_{d_n}$.
By the base chage of $\hat B(x) \to B(x)$ by $B(\xi) \to B(x)$,
we get a morphism 
$b(\xi):\hat B(\xi) \to B(\xi)$. 
The power series $f$ defines a closed subscheme $Z(f(\xi))$ of $B(\xi)$.
Then $\hat B(\xi)$ is quasi-smooth and the reduced part of 
$(b(\xi))^{-1}(Z(f(\xi)) \cup\cup_{i=1}^nZ(\xi_i))$ is a
quasi-normal crossing divisor:
$$
((b(\xi))^{-1}(Z(f(\xi)) \cup\cup_{i=1}^nZ(\xi_i)))_{red}
=\cup_{i=0}^{s+n}D_i,
$$
where $D_0$ and $D_{s+1}, \dots , D_{s+n}$ are proper transforms of
$Z(f(\xi))$ and $Z(\xi_1), \dots ,$ \linebreak
$Z(\xi_n)$ respectively and
$D_1, \dots , D_s$ correspond to 1-dimensional cones $r_1, \dots r_s$
of $F$ different from $\bold R_+e_i^*$ ($i = 1, \dots , n$).
Let $B(t) = Spec (k[[t]])$.
The morphism $B(x) \to B(t)$ (resp. $B(\xi ) \to B(t)$)
defined by $t=f(x_1, \dots , x_n)$ (resp. $t= f(\xi_1^{d_1}, \dots ,
\xi_n^{d_n})$) is denoted by $f(x)$ (resp. $f(\xi )$).  
$$
\CD
\hat B(\xi) @>>> \hat B(x) \\
@V{b(\xi)}VV @VV{b(x)}V \\
B(\xi) @>>> B(x) \\
@V{f(\xi)}VV @VV{f(x)}V \\
B(t) @>{=}>> B(t) \\
\endCD
$$
Let $\hat f(\xi ) = f(\xi)\circ b(\xi )$.
The multiplicities of $D_1, \dots , D_s$ in 
$b(\xi)^{-1}(Z(f(\xi)))= (\hat f (\xi))^{-1}(0)$ are computed as follows.
Let $L(x) = \bold Z^n$ and $L(\xi)$ be the lattice defined by
$\frac{1}{d_1}\bold Z \oplus \cdots \oplus \frac{1}{d_n}\bold Z$.
The dual lattice of $L(x)$ and $L(\xi)$ are denoted by
$L(x)^*$ and $L(\xi)^*$, respectively.  Then we have
$$
L(x) \subset L(\xi) \subset \bold R^n,
L(\xi)^* \subset L(x)^* \subset (\bold R^n)^*.
$$
The primitive generator of $r_i$ with respect to the dual lattice
$L(\xi)^*$ is denoted by $l_i$.  The multiplicity $m_i$ of $D_i$ 
in $(b(\xi))^{-1}(Z(f(\xi)))$ is equal to 
$$
m_i = \min \{ l_i (x) \mid x \in \Delta \}.
$$
\demo{Definition 1.1.1}
Let $p$ be the characteristic of $k$ and 
$d_1, \dots , d_n$ be natural numbers prime to $p$.
If there exists a simplicial refinement $F$ of $F(\Delta)$
with the property (1.1.1) such that $m = \lcm (m_i)$
is prime to $p$, the $n$-tuple of natural numbers $(d_1, \dots , d_n)$
is said to be tame with respect to $f$.
The number $m$ is called the exponent of $(d_1, \dots , d_n)$
with respect to $f$. (We fix the simplicial refinement $F$ of
$F(\Delta )$ once and for all.)
\enddemo

\heading
\S 1.2 Reduction of multiplicity
\endheading

We reduce the multiplicity of the special fiber
of $f(\xi)$ by the base change with respect to the morphism
$B(\tau) \to B(t)$.
We use the same notations of \S 1.1.  Let $m = \lcm (m_i)$
and $B(\tau ) = Spec (k[[\tau]])$, where $\tau^m = t$.  We
assume $(m,p) = 1$, i.e. $(d_1, \dots , d_n)$ is tame with respect
to $f$.  We construct a variety $\tilde B(\xi)$
over $B(\tau )$ and the following commutative diagram:
$$
\CD
\tilde  B(\xi) @>{\pi}>> \hat B(\xi)  \\
@V{\tilde f}VV @VV{\hat f}V \\
B(\tau ) @>>> B(t) \\
\endCD
$$
such that (1) $\tilde f^{-1}(0)$ is reduced and quasi-normal crossing
and (2) $\pi$ is a finite Galois covering with a Galois group 
$\mu_m$.

Let $F$ be a simplicial refinement of  $F(\Delta)$. 
First we construct a fan $\tilde F$ in $(\bold R^{n+1})^*$,
the suspension of $F$.  The polar dual $\Delta^*$ of the Newton
polyhedron $\Delta$ is defined by
$$
\Delta^* = \{x^* \in (\bold R^n)^* \mid
x^*(\Delta ) \geq 1\}.
$$
Since $F$ is a refinement of $F(\Delta )$, $F$ defines a decomposition
$\Delta^* = \coprod_{\sigma^* \in F}\sigma^* \cap \Delta^*$
of $\Delta^*$.  We define lower cone $L(\sigma^*)$,
upper cone $U(\sigma^*)$ of $\sigma^*$, 
and boundary cone $B(\sigma^*)$ as 
$$
\align
& L(\sigma^* ) =
\bold R_+((\sigma^* \cap \Delta^*)\times [0,\frac{1}{m}]),
U(\sigma^*) = (\sigma^* \times [0, \infty]-L(\sigma^*))^{-}, \\
& B(\sigma^*) = \bold R_+((\partial \Delta^*\cap
\sigma^*), \frac{1}{m}). \\
\endalign
$$
The base cone $B$ is defined by the closure of
$(\bold R^n)^*-\cup_{\sigma^* \in F}U(\sigma^*)$.
We define the suspension $\tilde F$ of $F$ as 
$$
\tilde F = \{U(\sigma^*)\}_{\sigma^* \in F} \cup
\{B(\sigma^*)\}_{\sigma^* \in F} \cup \{ B\}.
$$
It is easy to see that if $\sigma^*$ is
a simplicial cone, $U(\sigma^*)$ and $B(\sigma^* )$ are
also simplicial cones.  Even though if $F$ is a simplicial fan,
$B$ may not be simplicial cone.
Since the suspension $\tilde F(\Delta )$ of $F(\Delta )$ is 
the dual fan $F(\tilde \Delta )$ of the suspension $\tilde \Delta$,
where $\tilde \Delta$ is the convex hull of $0\times \bold R_{\geq m}$
and $\Delta \times 0$.  In general, $\tilde F$ is 
a refinement of $F(\tilde \Delta )$.
Now we consider an element
$f(\xi )-\tau^m=f(\xi_1^{d_1}, \dots ,\xi_n^{d_n})-\tau^m$ 
in $k[[\xi_1, \dots ,\xi_n,\tau ]]$.  
Let $B(\xi, \tau) = Spec(k[[\xi_1, \dots ,\xi_n,\tau]])$.
The fiber product of $X_{\tilde F} \to \bold A(\xi, \tau)$
and $B(\xi,\tau) \to \bold A(\xi, \tau)$ is denoted by 
$\bar B(\xi )$.
Since $\tilde \Delta$ is the Newton
polygon of $f(\xi )-\tau^m$, $f(\xi )-\tau^m$ 
defines a hypersurface $\tilde B(\xi) = Z(f(\xi)-\tau^m)$
in $\bar B(\xi)$.
The composite of $\tilde B(\xi)  \to \bar B(\xi )$,
$\bar B(\xi ) \to B(\xi, \tau)$ and
$B(\xi, \tau) \to B(\tau )$ is denoted by $\tilde f(\xi)$.
\proclaim{Proposition 1.2.1 (cf. [Dan])}
The variety $Z(f(\xi)-u^m)$ does not contain the point corresponding
to the base cone $B$ in $\bar B(\xi )$ and quasi-smooth.
The ( total ) fiber $(\tilde f (\xi))^{-1}(0)$ of $\tilde f$ at $0$ is
reduced and quasi-normal crossing.
\endproclaim

\heading
\S 2 Tame $l$-adic sheaves and convolution theorem.
\endheading
\heading
\S 2.1 Grothendieck groups of etale sheaves.
\endheading

In this section, we assume that $k$ is a finite field $\bold F_q$ and
$d$ be an integer such that $d \mid q-1$.  The Grothendieck group
of etale $\bold Q_l$-sheaves on $B(x)=Spec(k[[x]])$, $B^0(x)=Spec(k((x)))$ 
and $Spec (k)$ 
is denoted by $K(B(x))$,
$K(B^0(x))$ and $K(k)$, respectively. 
The geometric generic point of $B^0(x)$ is denoted by $\bar x$.
A sheaf $F$ on $B^0(x)$ corresponds to 
the continuous representation $(\rho_F, F_{\bar x})$ of 
$\Gal (\overline {k((x))}/k((x)))$.  
By Grothendieck's theorem, there exists an open subgroup $H$ in the inertia
group $I$ of $\Gal (\overline {k((x))}/k((x)))$
which acts on $F_{\bar x}$ unipotently.  
If there exists a subgroup $H$ such that the index $[I: H]$  
is prime to $p$, the sheaf $F$ is called tame.
The Grothendieck group of tame sheaves on $B^0(x)$ is denoted by $K(B^0(x))^t$.
For a tame sheaf $F$, $\min \{[I:H] \mid ([I:H],p) = 1, \text{ and
$\rho_F\mid_{H}$ is unipotent }\}$ is called the exponent of $F$.
The Grothendieck group of tame sheaves on $B^0(x)$ whose exponents divide $d$ 
is denoted by $K(B^0(x))^{t,d}$.  
Note that $K(B^0(x))^{t,1}$ is nothing
but the Grothendieck group of etale sheaves on $B^0(x)$ whose inertia
action is unipotent.  
Let us fix a generator $\gamma$ of the tame quotient 
$\hat\bold Z(1)'=\prod_{(p,l)=1}\bold Z_l(1)$ of $I$.
Then the logarithm $N= \log \rho_F (\gamma )$
of $\rho_F (\gamma )$ acts on $F_{\bar x}$ nilpotently.
The action of $N$ on $F_{\bar x}$ defines a monodromy filtration
$W_k(F_{\bar x})$ on $F_{\bar x}$ satisfying
$NW_k(F_{\bar x}) \subset W_{k-2}(F_{\bar x})$.
It is characterized by the property:
$$
\align
& \text{The $k$-th iteration 
$N^k: Gr_k^W(F_{\bar x}) \to Gr_{-k}^W(F_{\bar x})$
of the homomorphism} \\
& \text{ induced by $N$ is an isomorphism.} \\
\endalign
$$
If we define the primitive part $P_k$ as the kernel of the homomorphism
$N^{k+1}: Gr_k^W(F_{\bar x}) \to Gr_{-k-2}^W(F_{\bar x})$,
then $F$ is equal to 
$\sum_{k\geq 0}\sum_{0 \leq i \leq k}N^iP_k$ in $K(B^0(x))^{t,1}$.
Since $N_iP_k \simeq P_k \otimes \bold Q_l(i)$
 ($0 \leq i \leq k$) corresponds to an unramified representation
of $\Gal (\overline {k((x))}/k((x)))$,
we have the following lemma.
\proclaim{Lemma 2.1.1}
The group $K(B^0(x))^{t,1}$ is isomorphic to the group 
$K(k)$.
\endproclaim

  Let $\Pi_x^t = \Gal ((\cup_{(d,p)=1}k((x^{\frac{1}{d}})))\otimes \bar k
/k((x)))$
be the tame quotient of the absolute Galois group
of $k((x))$.  Then we have the exact sequence:
$$
1 \to \hat\bold Z(1)' \to \Pi^t_x  \to \hat \bold Z \to 1.
$$
The Kummer extension $k((\tau))$, where $\tau^d = x$, defines a character 
$\Pi^t_x \to \mu_d$ of $\Pi_x^t$.  The kernel of 
$\Pi_x^t \to \mu_d \times \hat\bold Z$ can be identified with 
$d \hat\bold Z (1)'$. We fix a generator $\gamma$ of $d\hat\bold Z(1)'$.
For a tame etale sheaf $F$ on $B^0(x)$ whose exponent is divisible by $d$,
the action of $d\hat\bold Z(1)'$ is unipotent.  Let $N$ be the logarithm
of $\rho_F(\gamma )$.  For any element $\bar h \in \mu_d\times \hat\bold Z$
take a lifting of $h \in \Pi^t_x$.  Then we have 
$h\gamma h^{-1} = \gamma^{c(h)}$, $c(h) \in (\hat\bold Z')^{\times}$
and as a consequence, we have $hN = c(h) Nh$.  If we define a filtration 
$W_i' = hW_i$, then $W_i'$ is stable under the action of $N$ and the 
following diagram is commutative.
$$
\CD
W_i/W_{i-1} @>{N}>> W_{i-2}/W_{i-3} \\
@V{c(h)\cdot h}VV @VV{h}V \\
W_i'/W_{i-1}' @>{N}>> W_{i-2}'/W_{i-3}' \\
\endCD
$$
By the characterization of $W_i$, we have $W_i' = W_i$.  Therefore $h$
induces a homomorphism on $Gr_i^W F_{\bar x}$.  Since $N$ acts trivially
on $Gr_i^WF_{\bar x}$, the action of $\Pi^t_x$ factors through 
$\mu_d \times \hat\bold Z$.  As a consequence, $K(B^0(x))^{t,d}$
is equal to the Grothendieck group of continuous representations of 
$\mu_d \times \hat\bold Z$.  Let $K(k, \mu_d)$ be the Grothendieck group
of etale sheaves on $Spec (k)$ with a $\mu_d$ action.
Therefore we have the following lemma.
\proclaim{Lemma 2.1.2}
The group $K(B^0(x))^{t,d}$ is isomorphic to the group 
$K(k, \mu_d)$. 
\endproclaim

We introduce a notion of an equivariant sheaf.  Let $f : X \to S$ 
be a scheme over $S$
and $G$ be a finite group acting on $X$ over $S$. 
The action of $G$ on $X$ is
denoted by $\sigma : G \to \Aut (X/S)$.  Let $\Cal F$ be a sheaf
on $X$.  A descent data for $\Cal F$ is a set of sheaf homomorphisms
$\phi_g : g^* \Cal F \to \Cal F$ indexed by elements of $G$ satisfying
the 1-cocycle condition $\phi_{gh} = \phi_h\circ h^*\phi_g $.  The descent
data is called effective if for any $x \in X$ and $g \in G$ such that
$g(x) = x$, the fiber of the descent data $\phi_g\mid_{x} :
(g^*\Cal F)_x= \Cal F_{g(x)} \to \Cal F_x$ is the identity map.
(Note that this terminology is different from the usual one.)  A sheaf
on $X$ with an effective descent data is called a $(G, \sigma )$-sheaf
on $X$. The Grothendieck group of $(G, \sigma )$-sheaves on $X$
is denoted by $K(X/S, (G, \sigma ))$.  Let $\Cal F$ be a $(G, \sigma )$-sheaf
on $X$.  Then the higher direct image sheaves $R^i f_* \Cal F$ is
a sheaf on $S$ with an action of $G$.
Suppose that there exists a
quotient scheme $Y=X/G$ of $X$ under the action of $G$. The structure
morphism $Y \to S$ is denoted by $g$ and the
natural map $X \to Y$ is denoted by $\pi$.  For an etale sheaf
$\Cal G$ on $Y$, the pull back $\pi^* \Cal G$ of $\Cal G$ by the morphism 
$\pi$ has a natural effective descent data.
It is known that
$(R^if_*(\pi^*\Cal G))^G \simeq R^ig_* \Cal G$.

\proclaim{Lemma 2.1.3}
The group $K(B^0(x))^{t,d}$ is isomorphic to the 
Grothendieck group 
$K(B^0(\tau), (\mu_d, \sigma))^{t,1}$ of 
unipotent $(\mu_d, \sigma )$-sheaves on $B^0(\tau)$.
\endproclaim
\demo{Remark 2.1.4}
This isomorphism depends on the choice of the uniformizer $x$ of
$k((x))$ which is always fixed in our context.
\enddemo

Let $F$ be an etale sheaf on $B(x)$.  The sheaf $F$ is said to be
tame if the restriction to $B^0(x)$ is tame.  For a tame sheaf $F$,
the exponent of $F$ is defined by that of the restriction to
$B^0(x)$.  The Grothendieck group of tame sheaves and tame sheaves whose
exponent divides $d$ are denoted by $K(B(x))^t$ and $K(B(x))^{t,d}$,
respectively.  For an etale sheaf $F$ on $B(x)$, the
generic geometric fiber and special geometric 
fiber are denoted by $F_{\bar x}$
and $F_{\bar 0}$ respectively.  The specialization map is
a $Gal(\bar k/ k)$-equivariant map 
$$
sp_F: F_{\bar 0} \to
F_{\bar x}^{I},
\tag{2.1.1}
$$ 
where $I$ is the inertia group of 
$\Gal (\overline{k((x))}/k((x)))$.  To give an etale sheaf on $B(x)$ is
equivalent to give a triple $(F_{\bar x}, F_{\bar 0}, sp_F )$,
where $F_{\bar x}$ and $F_{\bar 0}$ are 
a $\Gal (\overline{k((x))}/k((x)))$-module and a $\Gal (\bar k/ k)$-module
respectively and $sp_F$ is a $\Gal (\bar k/ k)$-equivariant homomorphism
(2.1.1).
By considering exact sequence
$$
0 \to j_!j^*F \to F \to i_*i^* F \to 0,
$$
we have $[F] = [j_!j^* F] + [i_*i^* F]$ in $K(B(x))$.  Since
$(j_!j^*F)_{\bar 0} = 0$ and $(i_*i^* F)_{\bar x} = 0$, we have
$$
\align
& K(B(x)) = K(B^0(x)) \oplus K(k),
K(B(x))^t = K(B^0(x))^t \oplus K(k), \\
& K(B(x))^{t,d} = K(B^0(x))^{t,d} \oplus K(k). \\
\endalign
$$

Now we consider several Grothendieck groups of etale
sheaves on $B(x) = $ \linebreak
$Spec(k[[x_1, \dots ,x_n]])$.  For a subset $I$ of
$[1,n]$, we define a closed subscheme 
$B_I(x) = Spec(k[[x_1,\dots , x_n]]/(x_i)_{i \in I})$ and locally closed 
subscheme $B_I^0(x) = B_I(x) - $ \linebreak
$\cup_{J\supsetneq I}B_J(x)$.  
Note that $B_{\emptyset}(x) = B(x)$.
A sheaf is said 
to be stratified by coordinate if the restriction to $B_I^0$ is locally
constant for all $I$.  An etale sheaf $F$ on $B(x)$ stratified by coordinate
is called tame if the monodromy of the restriction of $F$ to $B_I^0$
is tame.  The Grothendieck group of tame etale sheaves on $B(x)$ 
stratified by coordinate is denoted by $K_c(B(x))^t$.  
The tame fundamental
group of $B_I^0(x)$ is denoted by $\Pi_I^t$.  Then we have the following
 exact sequence.
$$
1 \to (\hat\bold Z(1)')^{I^c} \to \Pi_I^t \to \hat\bold Z \to 1 
$$
Let $I \subset [1,n]$ and $i \notin I$.  We set $J = I\cup {i}$ and fix
a generator $e_i$ of geometric monodromy along the divisor $x_i =0$.
Then we have
$1 \to \hat\bold Z(1)' e_i \to \Pi_I^t \to \Pi_J^t \to 1$.
The Grothendieck group of tame etale $\bold Q_l$-sheaves stratified
by coordinate 
whose exponent for $e_i$ divides $d_i$ is denoted by
$K_c(B(x))^{t, d_1, \dots ,d_n}$.
The generic geometric point of $B^0_I(x)$ is denoted by $\bar x_I$.
The specialization map with respect to $I$ and $J$ is 
a $\Pi_J^t$-equivariant map $sp_{I,J}: F_{\bar x_J} \to
F_{\bar x_I}^{e_i}$.  Let $K(\Pi_I^t)$ be the Grothendieck group
of continuous representations of $\Pi_I^t$.  By the same argument
as in the 1-variable case, we have
$$
K_c(B(x))^t = \oplus_{I \subset [1,n]}K(\Pi_I^t).
$$
If $F_1, \dots, F_n$ are tame sheaves on $B(x_1), \dots , B(x_n)$,
then the exterior product $F_1 \boxtimes \cdots \boxtimes F_n =
pr_1^* F_1 \otimes \cdots \otimes pr_n^* F_n$ is a tame sheaf on $B(x)$
stratified by coordinate.
\proclaim{Lemma 2.1.5}
By attaching $[F_1]\otimes \cdots \otimes [F_n]$ 
to $F_1 \boxtimes \cdots \boxtimes F_n$,
we get a homomorphism
$$
\align
& K(B(x_1))^t \otimes \cdots \otimes K(B(x_n))^t \to K_c(B(x))^t \\
& K(B(x_1))^{t,d_1} \otimes \cdots \otimes K(B(x_n))^{t,d_n} \to 
K_c(B(x))^{t,d_1, \dots ,d_n}. \\
\endalign
$$
\endproclaim
\demo{Proof}
It is enough to prove that if $0 \to F_1' \to F_1 \to F_1'' \to 0$
is an exact sequence of etale sheaves on $B(x_1)$,
$$
0 \to F_1' \boxtimes F_2 \boxtimes \cdots \boxtimes F_n \to
F_1 \boxtimes F_2 \boxtimes \cdots \boxtimes F_n \to
F_1'' \boxtimes F_2 \boxtimes \cdots \boxtimes F_n \to 0
$$
is an exact sequence of tame etale sheaves on $B(x)$.  It is clear by the
description of tame etale sheaves on $B(x)$ stratified by coordinate.
\enddemo
Let $B(\xi_i) \to B(x_i)$ be the $\mu_{d_i}$-covering defined 
in \S 1.1 and $G$ be the group $\mu_{d_1} \times \cdots \mu_{d_n}$ which
acts on $B(\xi )$. The action is denoted be $\sigma_\xi$.
\proclaim{Proposition 2.1.6}
\roster\item
$
K_c(B(\xi), (G, \sigma_\xi))^{t,1} \simeq
\oplus_{I \subset [1,n]} K(B_I^0(\xi), (G_I, \sigma_I))^{t,1},
$
where $G_I$ is the image of $G$ in $\Aut (B_I^0(\xi))$ and
$\sigma_I$ be the action of $G_I$ on $B^0_I(\xi)$, i.e.
$G_I = \prod_{i \notin I}\mu_{d_i}$.
\item
$K_c(B_I^0(\xi), (G_I, \sigma_I))^{t,1} \simeq K(k, G_I)$.
\endroster
\endproclaim

\heading
\S 2.2 Vanishing cycle functors and convolution theorem
\endheading

  Let $I$ be a subset of $[1,n]$.  We use the same notations for
$B_I(x)$, $B(\xi)$, $G$ and $\sigma_\xi$ as \S 2.1.
Let $[F_1]\otimes \cdots \otimes [F_n]$ be an element of 
$K(B(x_1))^t\otimes \cdots \otimes K(B(x_n))^t$.
We define a homomorphism
$\phi_I : K(B(x_1))^{t,d_1} \otimes \cdots \otimes K(B(x_n))^{t, d_n}
\to K(B^0(t))$ as
$$
\phi_I([F_1] \otimes \cdots \otimes [F_n])
 = [j_t^*(\bold R f(x)_*{i_I}_*i_I^*(F_1 \boxtimes \cdots \boxtimes F_n))],
$$
where $j_t : B^0(t) \to B(t)$ and $i_I: B_I(x) \to B(x)$ be the
natural inclusions. 
For an element 
$[\tilde F_1]\otimes \cdots \otimes [\tilde F_n]\in
K(B(\xi_1), (\mu_{d_1}, \sigma_1))^{t} \otimes \cdots \otimes
K(B(\xi_n), (\mu_{d_n}, \sigma_n))^{t}$,
we define the equivariant version
$$
\tilde \phi_I :K(B(\xi_1), (\mu_{d_1}, \sigma_1))^{t,1} \otimes \cdots \otimes
K(B(\xi_n), (\mu_{d_n}, \sigma_n))^{t,1}\to K(B^0(t),G),
$$
of $\phi_I$ by
$$
\tilde \phi_I ([\tilde F_1]\otimes \cdots \otimes [\tilde F_n]) =
[j_t^* \bold R f(\xi)_*\tilde i_{I*}\tilde i_I^*
(\tilde F_1 \boxtimes \cdots \boxtimes \tilde F_n)],
$$
where $\tilde i_I: B_I(\xi) \to B(\xi)$ is the natural inclusion.
Then we have the following commutative diagram.
$$
\CD
K(B(x_1))^{t,d_1} \otimes \cdots \otimes K(B(x_n))^{t,d_n}
@>{\phi_I}>>
K(B^0(t)) \\
@VVV @AA{\text{ $G$-invariant }}A\\
K(B(\xi_1), (\mu_{d_1}, \sigma_1))^{t,1} \otimes \cdots \otimes
K(B(\xi_n), (\mu_{d_n}, \sigma_n))^{t,1} @>>{\tilde \phi_I}>  K(B^0(t),G) \\
\endCD
\tag{2.2.1}
$$
First we prove the following proposition.
\proclaim{Proposition 2.2.1}
The image of $\phi_I$ is contained in $K(B^0(t))^{t,m}$, where $m$
is the exponent of $(d_1, \dots , d_n)$ with respect to $f$.
\endproclaim
\demo{Proof}
We consider the following commutative diagram:
$$
\CD
\tilde B(\xi) @>{\nu}>> \hat B(\xi) @>{b}>> B(\xi) \\
@V{\tilde f(\xi)}VV @VV{\hat f(\xi)}V \\
B(\tau) @>>{\nu '}> B(t). \\
\endCD
$$
For tame etale sheaves 
$F_1, \dots , F_n$ on $B(x_1), \dots , B(x_n)$ whose exponent divides
$d_1, \dots ,d_n$, the sheaf $\pi_1^*F_1 \boxtimes \cdots \boxtimes \pi_n^*F_n$
are $(G, \sigma_\xi)$-sheaf on $B(\xi)$ with unipotent monodromy.
Therefore there exists a filtration $F^i$ of $(G, \sigma_\xi)$-sheaves
on $\pi_1^*F_1 \boxtimes \cdots \boxtimes \pi_n^*F_n$ such that $F^i / F^{i+1}$
has trivial monodromy.  Therefore the subquotient 
$(b \circ \nu)^*F^i/ (b\circ\nu)^*F^{i+1}$ of 
$(b\circ\nu)^*(\pi_1^*F_1 \boxtimes \cdots \boxtimes \pi_n^*F_n)$ 
has also trivial monodromy.
Using weight spectral sequence for $\tilde f (\xi)$, 
$j_{\tau}^*R^i\tilde f(\xi)_*
((b\circ\nu)^*(\pi_1^*F_1 \boxtimes \cdots \boxtimes \pi_n^*F_n))$
has unipotent monodromy.  Since 
$$
j_{\tau}^*R^i\tilde f(\xi )_*
((b\circ\nu)^*(\pi_1^*F_1 \boxtimes \cdots \boxtimes \pi_n^*F_n))
\simeq
(\nu ')^*
j_{t}^*R^i f(\xi)_*
b^*(\pi_1^*F_1 \boxtimes \cdots \boxtimes \pi_n^*F_n),
$$
and the commutative diagram (2.2.1), we have
$$
[j_{t}^*R^i f(x)_*
(F_1 \boxtimes \cdots \boxtimes F_n)]=
[j_{t}^*R^i f(\xi)_*
b^*(\pi_1^*F_1 \boxtimes \cdots \boxtimes \pi_n^*F_n)]^G \in 
K(B^0(t))^{t,m}
$$
\enddemo
\demo{Definition 2.2.2}
Put $F = [F_1] \otimes \cdots \otimes [F_n]$.
\roster\item
$$
\Phi_f(F) =
\phi_{\emptyset}(F) - 
[(F_1 \boxtimes \cdots \boxtimes F_n)_{\bar 0}]
\in K(k,\mu_m).
$$
Here we use the identification, $\phi_{\emptyset}(F)
\in K(B(t)^0)^{t,m} \simeq K(k, \mu_m)$ and
$[(F_1 \boxtimes \cdots \boxtimes F_n)_{\bar 0}] \in K(k) = K(k, {1})
\subset K(k, \mu_m)$.
\item
$$
\align
& \Phi ([F_i]) = [F_{i,\bar x_i}] - [F_{i,\bar 0}] 
\in K(k, \mu_{d_i})\text{ and }
\\
& \tilde\chi_I(F) = \otimes_{i \notin I } [F_{i,\bar x_i}] \otimes
\otimes_{i \in I}\Phi ([F_i]) \in K(k, G), \\
& \chi_I(F) = \otimes_{i \notin I } [F_{i,\bar x_i}] \otimes
\otimes_{i \in I}\Phi ([F_i])^{\mu_{d_i}} \in K(k, G_I), \\
\endalign
$$
\item
$$
\align
\tilde\Phi_I  (F) =& 
\tilde\phi_I([\pi_1^*F_1]\otimes\cdots\otimes [\pi_n^* F_n])
-[(F_1 \boxtimes \cdots \boxtimes F_n)_{\bar 0}] \\
& \in K(B^0(t), (G_I, \sigma_I))^{t,m} \simeq 
K(k, G_I \times \mu_m), \\
\endalign
$$
and
$\tilde \Phi_I = \tilde \Phi_I([\bold Q_l]\otimes\cdots\otimes [\bold Q_l])$.
\endroster
\enddemo
\proclaim{Theorem 2.2.3 ($l$-adic convolution theorem)}
Let $F = [F_1] \otimes \cdots \otimes [F_n]$ in
$K(B(x_1))^{t,d_1} \otimes \cdots \otimes K(B(x_n))^{t,d_n}$.
Under the notation defined as above,
$$
\Phi_f(F) = \sum_{I \subset [1,n]}
(-1)^{\# I}(\tilde \Phi_I \otimes \chi_I (F))^{G_I}
\in K(\mu_m, k).
$$
\endproclaim
\demo{Proof}
Since the action of $G$ on $\tilde\Phi_I$ factors through $G_I$,
it is enough to show that
$$
\tilde \Phi_f(F) = 
\sum_{I \subset [1, n]}(-1)^{\# I}
(\tilde \Phi_I \otimes \tilde\chi_I(F)) \in K(G\times \mu_m, k),
$$
by the commutative diagram (2.2.1).
Let $\tilde F = [\pi_1^*F_1 \boxtimes \cdots \boxtimes \pi_n^*F_n]$.
Using decomposition 
$
\tilde F = \sum_{I \subset [1, n]}
\tilde j_{I!}(\tilde F\mid_{B_I^0(\xi)}),
$
where $\tilde j_I : B_I^0(\xi) \to B(\xi)$ is the natural inclusion.
Since $ \tilde F\mid_{B_I^0(\xi)}
 \in K(B_I^0(\xi), (G_I, \sigma_I))^{t,1}$,
we have
$$
\tilde j_{I !}( \tilde F\mid_{B_I^0(\xi)}) =
[\tilde j_{I !} \bold Q_l] \otimes
\otimes _{i \notin I}[F_{i, \bar x_i}]\otimes \otimes_{i \in I}
[F_{i, \bar 0}]
\in K_c(B(\xi), (G, \sigma_\xi))^{t,1}.
$$
Here the $(G, \sigma_\xi)$-structure on the right hand side is given by the
diagonal action of $(G_I, \sigma_I)$ on $\otimes_{i\notin I}[F_{i, \bar x_i}]$
and $[\tilde j_{I !}\bold Q_l]$.  Using equality $[\tilde j_{I !} \bold Q_l] =$
\linebreak
$\sum_{I \subset K}(-1)^{\# (K-I)}[\tilde i_{K !}\bold Q_l]$, where
$\tilde i_K : B_K(\xi) \to B(\xi)$ is the natural inclusion, we have
$$
\align
F_B =& \sum_{K \supset I}(-1)^{\# (K-I)}
[\tilde i_{K!} \bold Q_l] \otimes \otimes_{i \notin I} [F_{i,\bar x_i}]
\otimes \otimes_{i \in I}[F_{i, \bar 0}] \\
= &
 \sum_{K \subset [1,n]}(-1)^{\# K}
[\tilde i_{K!} \bold Q_l] \otimes 
(\sum_{K \supset I}(-1)^{\# I}
\otimes_{i \notin I} [F_{i,\bar x_i}]
\otimes \otimes_{i \in I}[F_{i, \bar 0}] ). \\
\endalign
$$
Since 
$$
\align
& \sum_{K \supset I}(-1)^{\# I}
\otimes_{i \notin I} [F_{i,\bar x_i}]
\otimes \otimes_{i \in I}[F_{i, \bar 0}] \\
= &
\otimes_{i \notin K}  [F_{i,\bar x_i}]
\otimes 
(\sum_{K \supset I}(-1)^{\# I}
\otimes_{i \notin K-I}  [F_{i,\bar x_i}] 
\otimes \otimes_{i \in I} [F_{i, \bar 0})] \\
= &
\otimes_{i \notin K} [F_{i,\bar x_i}]
\otimes 
\otimes_{i \in K} \Phi([F_i]) \\
= &
\tilde\chi_K(F), \\
\endalign
$$
we have
$$
\align
\tilde\Phi_f(F) & = \sum_{K \subset [1, n]}
(-1)^{\# K}\tilde\Phi_f([\tilde i_{K !}\bold Q_l]) 
\otimes \tilde\chi_K(F) \\
& = \sum_{K \subset [1, n]}
(-1)^{\# K}\tilde\Phi_K 
\otimes \tilde\chi_K(F). \\
\endalign
$$
This completes the proof of the theorem.
\enddemo

We apply Theorem 2.2.3 to composite singularities.
Let $m_1, \dots, m_n$ be positive integers, 
$y_i = (y_{i1}, \dots , y_{im_i})$ sets of coordinates.  We  define 
$B(y_i) = $ \linebreak
$Spec(k[[ y_{i1}, \dots , y_{im_i}]])$.
Let $g_i = g_i(y_i)$ be a formal power series of $\{y_{ij}\}_j$ with
no constant term.  Then $g_i$ defines a homomorphism $B(y_i) \to B(x_i)$
by $x_i = g_i(y_i)$.  We assume that $g_i$ is smooth on 
$B(y_i) - g_i^{-1}(0)$ and the exponent $d_i$ of $g_i$ divides
$\# k^{\times}$.
The fiber product of $g_i$ defines a morphism from
$B(y) = \prod_{i=1}^nB(y_i)$ to $B(x) = \prod_{i=1}^nB(x_i)$ and it
is denoted by $g$.  
We assume that the exponent $m$ of $(d_1, \dots ,d_n)$ with respect to
$f$ divides $\# k^{\times}$.
$$
B(y)=\prod_{i=1}^n B(y_i) \overset{g =\prod_{i=1}^ng_i} \to\longrightarrow
B(x)=\prod_{i=1}^n B(x_i) \overset {f} \to\to B(t)
$$
\proclaim{Lemma 2.2.4}
The composite $f\circ g : B(y) \to B(x) \to B(t)$ is smooth over
$B^0(t)$.
\endproclaim
Let $\chi_{g,K}(\bold Q_l) = \chi_K([\bold Rg_{1*} \bold Q_l]\otimes
\cdots \otimes [\bold R g_{n*} \bold Q_l])$ and
$$
\Phi_{f\circ g}(\bold Q_l)
=\Phi_f( [\bold Rg_{1*}\bold Q_l]\otimes \cdots\otimes
[\bold R g_{n*}\bold Q_l])
=[j_t^*\bold R(f\circ g)_* \bold Q_l]- [\bold Q_l].
$$
Then the following theorem is a direct consequence of Theorem 2.2.3.
\proclaim{Theorem 2.2.5 (l-adic convolution theorem for composite 
singularities)}
$$
\Phi_{f\circ g}(\bold Q_l) = 
\sum_{I \subset [1,n]}(-1)^{\# I} 
(\tilde\Phi_I \otimes \chi_{g, I}(\bold Q_l))^{G_I}.
$$
\endproclaim
\demo{Remark 2.2.6}  If we apply this theorem for the function 
$f(x_1, x_2) = x_1 + x_2$, it is nothing but the $l$-adic version of
Thom-Sebastiani theorem. (See also Corollary 3.6.2.)  
It is proved in [L-D] assuming the existence
of resolution of singularities.  Here we only assume the tameness of
the singularity for $g_1=0$ and $g_2=0$. 
\enddemo

The rest of this paper is devoted to prove the Hodge analog of
Theorem 2.2.5.

\heading
\S 3 Mixed Hodge structure of composite singularities
\endheading

\heading
\S 3.1 Local space
\endheading

  Let $N$ be an analytic space and $X$ and $Y$ be closed subspaces
such that $X$ is a projective algebraic variety.  The triple $(N, X, Y)$
is called a local space if the following conditions hold.
\roster
\item
Any irreducible component of $Y$ is not contained in $X$.
\item 
$N-(X\cup Y)$ is a smooth variety.
\item
There exists a proper bimeromorphic map $b: \tilde N \to N$ such
that $b$ is an isomorphism on $b^{-1}(N-(X \cup Y))$ and
$b^{-1}(X \cup Y)$ is a normal crossing divisor of $\tilde N$.

The irreducible decomposition of the proper transform $\tilde Y$ of 
$Y$ is written by $\tilde Y = \cup_{i \in I}\tilde Y_i$.
For a subset $J$ of $I$,  
we define $\tilde Y_J = \cap_{j \in J} \tilde Y_j$. 
\item 
For any subset $J$ of $I$,
$\tilde Y_J$ retracts to $\tilde Y_J \cap \tilde X$,
where $\tilde X = b^{-1}(X)$.
(In particular, $\tilde N$ retracts to the subvariety $\tilde X$.)
\endroster
Note that for a triple with the properties (1)-(3), there exists
a sufficiently small neighborhood $N'$ of $X$ such that
$(N', X, Y\cap N'')$ is a local space.  We introduce an equivalence 
relation on the set of local spaces by $(N, X, Y) \sim (N', X, Y\cap N')$,
where $N'$ is a sufficiently small neighborhood of $X$.
From now on, equivalent local spaces are always identified.

\demo{Example 3.1.1}
Let $B(x) = \{ (x_1, \dots ,x_n)\mid \mid x_i\mid < \epsilon \}$
and $f(x_1, \dots , x_n)$ be a holomorphic function on $B(x)$
to $B(t) =\{t \mid \mid t\mid <\epsilon \}$.  Then for a sufficiently small
$\epsilon$, $(B(x), 0, f^{-1}(0))$ is a local space.
\enddemo
In this section, we introduce a mixed Hodge structure on the cohomology
$H^i(N- (X \cup Y), \bold Q)$.  To introduce a Hodge structure, we
fix a resolution of singularities $\tilde N \to N$ of $(3)$.  Let 
$\tilde Y =\cup_{i\in I}\tilde Y_i$ and $\tilde X = \cup \tilde X_i$
be the irreducible decompositions of the proper transform $\tilde Y$
of $Y$ and that of $b^{-1}(X)$, respectively.  
Then $(\tilde N, \tilde X, \tilde Y)$ is also a local space.
First we construct a mixed Hodge structure on 
$H^i(\tilde N -(\tilde X \cup \tilde Y), \bold Q)
=H^i(N -(X \cup Y), \bold Q)$.  Later we prove
that it is independent of the choice of the resolution of singularities.
The local space $(N,X,Y)$ is called normal crossing if 
$N$ is smooth and $X \cup Y$
is a normal crossing divisor in $N$.
We assume that $(N, X, Y)$ is a normal crossing local space and
$X = \cup_{i=1}^q X_i$ and $Y = \cup_{i=1}^pY_i$ are the 
irreducible decomposition.  We can define the logarithmic de Rham 
complex $\Omega^{\bullet}_N(\log D)$ as in [Del].
By the result of [Del], there exist quasi-isomorphisms
$\bold R j_*\bold C  \to j_*\Omega^{\bullet}_{N-D}$
and $\Omega_N(\log D ) \to j_*\Omega_{N-D}^{\bullet}$,
where $j: N-D \to N$ is the natural inclusion.
Let $i_I : X_I \to X$ be the natural inclusion.
If $I \subset J$, then there exists a natural morphism 
$(i_I)_* \bold Rj_* \bold Q \to (i_J)_* \bold Rj_* \bold Q$.

  Let $\Cal N$ be an analytic manifold and $D = \cup_iD_i$ be
a simple normal crossing divisor.  For a subset $J$ of $I$, 
we put $D_J = \cap_{j \in J}D_j$ and the ideal sheaf of $D_J$ is
denoted by $\Cal I_J$.  The sheaf of logarithmic differential
forms $\Omega^i_{\Cal N}(\log D)$ is a locally free sheaf on $\Cal N$
and $\Cal I_{D_J} \otimes \Omega^i_{\Cal N}(\log D)$ is identified
with a subsheaf of $\Omega_{\Cal N}^i(\log D)$.  If $z_i$
is a local equation of $D_j$ and $\omega \in \Omega_{\Cal N}^i(\log D)$,
then we have $d(z_ \omega )= z_i(\frac{dz_j}{z_j}\wedge \omega)+
z_jd\omega$, therefore
$\Cal I_{D_J} \otimes \Omega^{\bullet}_{\Cal N}(\log D)$ is a 
sub-complex of $\Omega_{\Cal N}^{\bullet}(\log D)$.
Therefore the quotient sheaves $\Omega^i_{\Cal N}(\log D) \otimes
\Cal O_{D_J}$ forms a complex of sheaves
$\Omega^{\bullet}_{\Cal N}(\log D) \otimes \Cal O_{D_J}$.
\proclaim{Lemma 3.1.2}
Let $j: \Cal N -D \to \Cal N$ be the natural inclusion. Then there
exists a quasi-isomorphism
$$
\bold R j_* \bold C\mid_{X_J} \simeq
\Omega_{\Cal N}^{\bullet}(\log D) \otimes \Cal O_{X_J}.
$$
\endproclaim
\demo{Proof}
First consider the following diagram where all the morphisms are filtered
quasi isomorphisms;
$$
(\bold R j_* \bold C, \tau_{\bullet})
\to
(\bold R j_* \Omega_{\Cal N- D}, \tau_{\bullet})
\longleftarrow
(\Omega_{\Cal N}(\log D), \tau_{\bullet})
\to 
(\Omega_{\Cal N}(\log D), W_{\bullet})
\tag{3.1.1}
$$
Let $W_k(\Omega_{\Cal N}(\log D)\otimes \Cal O_{D_J}) = 
\Im (W_k(\Omega_{\Cal N}(\log D)) \to
\Omega_{\Cal N}(\log D)\otimes \Cal O_{D_J})$.
Then we have the following lemma.
\enddemo
\proclaim{Lemma 3.1.3}
\roster
\item
The Poincare residue induces the isomorphism
$$
Gr_k^W(\Omega_{\Cal N}^i(\log D) \otimes \Cal O_{D_J})
\simeq \oplus_{\# K = k}\Omega^{i-k}_{D_{J_\cup K}}(-k)
$$
\item
The natural morphism $\pi :(\Omega_{\Cal N}^{\bullet}(\log D)
\mid_{D_J}, W_{\bullet})
\to (\Omega_{\Cal N}^{\bullet}(\log D)\otimes \Cal O_{D_{J}}, W_{\bullet})$
is a filtered quasi isomorphism
\endroster
\endproclaim
\demo{Proof}
(1) is well known and easy to check.

(2) The associated graded morphism of $\pi$ is equal to
$$
\oplus_{\# K = k}\Omega^{\bullet}_{D_K}\mid_{D_{J_\cup K}}[-k]
\to
\oplus_{\# K = k}\Omega^{\bullet}_{D_{J_\cup K}}[-k]
$$
and they are quasi-isomorphic to $\bold C_{D_{J \cup K}}$.  This proves the
lemma 3.1.3.
\enddemo
By (3.1.1) and Lemma 3.1.3 (2), we have Lemma 3.1.2.

Now we return to the situation of a normal crossing local space $(N, X,Y)$.
We define the following complexes $K_{\bold Q}=K_{\bold Q}(N,X,Y)$
and $K_{\bold C}=K_{\bold C}(N,X,Y)$ on $X$.
$$
K_{\bold Q} : \oplus_{\# I=1}(i_I)_*\bold R j_* \bold Q\mid_{X_I} \to
\oplus_{\# I=2}(i_I)_*\bold R j_* \bold Q\mid_{X_I} \to \cdots
$$
and
$$
K_{\bold C} : \oplus_{\# I=1}
(i_I)_*\Omega^{\bullet}_{\Cal N}(\log (D))\otimes\Cal O_{X_I} \to 
\oplus_{\# I=2}
(i_I)_*\Omega^{\bullet}_{\Cal N}(\log (D))\otimes\Cal O_{X_I} \to \cdots
$$
The following lemma is a direct consequence of Lemma 3.1.2.
\proclaim{Lemma 3.1.4}
There exists a quasi isomorphism of complexes of sheaves:
$$
\bold R j_* \bold Q \mid_X \to K_{\bold Q}, K_{\bold Q}\otimes \bold C \to K_C.
$$
\endproclaim
\demo{Definition 3.1.5}
We introduce a weight filtration $W$ on $K_{\bold Q}$, $K_{\bold C}$,
and a Hodge filtration $F$ on $K_{\bold C}$ as follows.
$$
\align
W_k K_{\bold Q}: & \oplus_{\# I=1}\tau_k\bold R {j}_* \bold Q
\mid_{X_I}
\to \oplus_{\# I=2}\tau_{k+1}\bold R {j}_* \bold Q\mid_{X_I}
 \\
& \to \cdots  \to \oplus_{\# I=j}\tau_{k+j-1}\bold R {j_I}_* \bold Q
\to \cdots \\
\endalign
$$
$$
\align
W_k K_{\bold C}: & \oplus_{\# I=1}W_k\Omega_{N}^{\bullet}
(\log (D))\otimes \Cal O_{X_I} 
\to \oplus_{\# I=2}W_{k+1}\Omega_{N}^{\bullet}(\log D)\otimes \Cal O_{X_I} \\
& \to \cdots  \to 
\oplus_{\# I=j}W_{k+j-1}\Omega_{N}^{\bullet}(\log D)\otimes \Cal O_{X_I} 
\to \cdots \\
\endalign
$$
$$
\align
F^p K_{\bold C}: & \oplus_{\# I=1}\sigma^p\Omega_{D}^{\bullet}
(\log (D))\otimes \Cal O_{X_I} \to 
\oplus_{\# I=2}\sigma^p\Omega_{N}^{\bullet}(\log D)\otimes \Cal O_{X_I} \\
& \to \cdots  \to 
\oplus_{\# I=j}\sigma^p\Omega_{N}^{\bullet}(\log D)\otimes \Cal O_{X_I} 
\to \cdots , \\
\endalign
$$
where $\tau_{\bullet}$ and $\sigma^{\bullet}$ are the canonical and
the stupid filtrations, respectively.
\enddemo
These filtrations on the complex of sheaves induces the weight filtration
$W$ on $\bold R\Gamma (X, K_{\bold Q})$ and 
$\bold R\Gamma(K, K_{\bold C})$ and the Hodge filtration 
$F$ on $\bold R \Gamma (X, K_{\bold C})$ up to filtered 
quasi-isomorphism.  
\demo{Definition 3.1.6}
The pair 
$
((L_{\bold Q}, W),(L_{\bold C}, W, F))$
of filtered complex of sheaves and bi-filtered complex of sheaves on X 
is a cohomological mixed Hodge complex (=CMHC for short [Del]), if
the following conditions hold.
\roster
\item
$(L_{\bold Q}, W)\otimes \bold C$
and 
$(L_{\bold C}, W)$
are filtered quasi-isomorphism.
\item
The spectral sequence for the filtered complex of sheaves
\linebreak
$(\bold R \Gamma (X, Gr_k^W(L_{\bold C})), Gr_k^W(F))$
on $X$ degenerates at $E_1$-term.
\item
The filtration $Gr_k^W(F)$
defines a Hodge structure of weight $m+k$ on 
\linebreak
$H^m(X, Gr_k^W(L_{\bold C}))$.
\endroster
\enddemo
\proclaim{Theorem 3.1.7}
The pair $((K_{\bold Q}, W), (K_{\bold C}, W,F))$ defines a 
CMHC on X.
\endproclaim
\demo{Proof}
It is sufficient to prove that $Gr_k^W(K_{\bold C})$
is a cohomological Hodge complex of weight $k$, where
$Gr_k^W(K_{\bold C})$ is given by the double complex:
$$
\oplus_{\# I=k,\# J =1}\Omega_{Y_I\cap X_J}^{\bullet}[-k](-k)  
\overset{0}\to\to
\oplus_{\# I=k+1,\# J =2}\Omega_{Y_I\cap X_J}^{\bullet}[-k-1](-k-1)  
\overset{0}\to\to 
$$
Moreover, the associated graded module $Gr_F^pGr_k^W(K_{\bold C})$
for the induced filtration
$Gr_k^W(F^p)$ is equal to
$$
\oplus_{\# I=k,\# J =1}\Omega_{Y_I\cap X_J}^{p-k}[-p](-k)  
\overset{0}\to\to
\oplus_{\# I=k+1,\# J =2}\Omega_{Y_I\cap X_J}^{p-k-1}[-p](-k-1)  
\overset{0}\to\to 
$$
Therefore the $m$-th hypercohomology of 
$Gr_F^pGr_k^W(K_{\bold C})$
is equal to
$$
\oplus_{j\geq 1}\oplus_{\# I =k+j-1, \# J =j}
H^{m-p-j+1}(\Omega_{Y_I\cap X_J}^{p-k-j+1})(-k-j+1).
$$
Therefore by the classical Hodge theory the following spectral
sequence degenerates at the $E_1$-term.
$$
E_1^{p,q}= \bold H^{p+q}(X, Gr_F^pGr_k^W(K_{\bold C}))
\to E_{\infty}^{p+q}=
\bold H^{p+q}(X, Gr_k^W(K_{\bold C})),
$$
and the filtration defined by this spectral sequence defines
a pure Hodge structure of weight $m+k$ on 
$\bold H^m(X,Gr_k^W(K_{\bold C}))$.
\enddemo

We prove the independence of the choice of the resolution
of singularity.  We use the same argument as [Del].  Let $(N_1, X_1, Y_1)$
and $(N_2, X_2, Y_2)$ be normal crossing local spaces,
$(N, X, Y)$  a local space and $f_i: N_i \to N$ be proper morphisms
for $i=1, 2$ such that $f^{-1}(X\cup Y)= X_i\cup Y_i$ and 
$f_i \mid_{N_i - (X_i\cup Y_i)}:
N_i -(X_i\cup Y_i) \to N -(X\cup Y)$ are isomorphic 
for $i=1, 2$.  Then by taking
a fiber product of $f_1, f_2$ and resolving it, there exists the
third normal crossing local space $(N_3, X_3, Y_3)$ and proper
morphism $g_i : (N_3, X_3, Y_3) \to (N_i, X_i, Y_i)$ such that 
$g_i^{-1}(X_i\cup Y_i) = X_3\cup Y_3$ and the restrictions of 
$g_i$ to $N_3-(X_3\cup Y_3)$ are isomorphisms.  
$$
\CD
(N_3, X_3, Y_3) @>>> (N_1, X_1, Y_1) \\
@VVV @VVV \\
(N_2, X_2, Y_2) @>>> (N, X,Y) \\
\endCD
$$
The morphism $g_i$
induces a homomorphism from the mixed Hodge structure of
$(N_i, X_i, Y_i)$ to that of $(N_3, X_3, Y_3)$ which is bijective.
Since the category of mixed Hodge structures is an abelian category,
this gives rise to the isomorphism of mixed Hodge structures.
In the same way, we can prove the following functoriality
of the mixed Hodge structures of local spaces.
\proclaim{Proposition 3.1.8}
Let $(N_1, X_1, Y_1)$ and $(N_2, X_2, Y_2)$ be local spaces and 
$f: N_1 \to N_2$ be a morphism such that $f(N_1 -(X_1 \cup Y_1))
\subset N_2- (X_2 \cup Y_2)$ and $f\mid_{N_1-(X_1\cup Y_1)}$
is an immersion.
Then $f$ induces a homomorphism of mixed Hodge structures of
local spaces.
\endproclaim
\demo{Proof}
We reduce the proof to the case where both $(N_1, X_1, Y_1)$ and
$(N_2, X_2, Y_2)$ are normal crossing local spaces.
In this case, the proof is exactly the same as [Del].
\enddemo
\demo{Remark 3.1.9}
Let $(N, X, Y)$ be a normal crossing local space and $X_1$ 
a closed subvariety of $X$,
$N_1$ a neighborhood of $X_1$,
and $Y_1$ the closure of 
$(X\cup Y) \cap N_1 -X_1$.  Then $(N_1, X_1, Y_1)$ 
is a local space and the open immersion $N_1 \subset N$
satisfies the condition (1).
\enddemo
\proclaim{Corollary 3.1.10}
Let $(N, X, Y)$ be a local space and $W$ be a subvariety of $X$.  Then
$H^i_c(W, \bold R j_*\bold Q)$ has a natural mixed Hodge structure.
\endproclaim
\demo{Proof}
Let $Z_1$ be the closure of $W$ in $X$ and put $Z_2 = Z_1 -W$.
The open immersion from $W$ to $Z_1$ is denoted by $k : W \to Z_1$.
Since 
$$
H^i_c(W, \bold R j_*\bold Q ) \simeq 
H^i(Z_1, k_!(\bold R j_* \bold Q \mid_W))
$$
and 
$$
k_!(\bold R j_*\bold Q\mid_W) \to \bold R j_*\mid_{Z_2}
\overset{rest}\to\to \bold R j_*\mid_{Z_1}
\overset{+1}\to\to
$$
is a distinguished triangle, and
by the construction of cone of CMHC given later in \S 3.4,
it is enough to prove that
$\bold R j_*\bold Q\mid_{Z_1}$ and
$\bold R j_*\bold Q\mid_{Z_2}$ are cohomological mixed Hodge complex
and the morphism $rest$ is a morphism of mixed Hodge complex.
For the first statement, we may assume that $W$ is the closed
subvariety of $X$.  By modifying, there exists a modification
of local spaces $\pi:(\tilde N, \tilde X, \tilde Y) \to (N, X, Y)$
such that $(\tilde N, \tilde X, \tilde Y)$ is normal crossing
and $\tilde W = \pi^{-1}(W)$ is a divisor of $\tilde N$ contained
in $\tilde X$.  By using the construction of Remark 3.1.9, we can see that
$\bold R\tilde j_*\bold Q \mid _{\tilde W}$ is a cohomological
mixed Hodge complex on $\tilde W$.
Therefore, 
$$
\bold R\Gamma (W, \bold R j_*\bold Q \mid_W)
\simeq \bold R\Gamma (W, \bold R \pi_*(\bold R\tilde j_* \bold Q
\mid_{\tilde W}))
\simeq \bold R \Gamma (\tilde W, \bold R \tilde j_* \bold Q \mid_{\tilde W})
$$
is a mixed Hodge complex.  We can prove that the morphism $rest$
is a morphism of CMHC's in the same way.
\enddemo

\heading
\S 3.2 Mixed Hodge structures associated to hypersurface
singularities
\endheading

In this section, we investigate mixed Hodge structures of germs of
hypersurface singularities.

Let $B(y) = \{(y_1, \dots , y_n) \mid \mid y_i \mid < \epsilon \}$
and $f = f(y_1, \dots , y_n)$ be a holomorphic function on $B(y)$
such that $f(0) =0$.  The function $f$ defines a holomorphic
map from $B(y)$ to $B(t) =\{ t \in \bold C\mid \mid t \mid < \epsilon\}$
and put $Z(f) = \{ y \in B(y) \mid f(y) =0\}$.  We assume that 
$f$ is smooth outside of $Z(f)$.
Let $m$ be a positive integer.  The variety $B_m(y)$ is defined by
the following fiber product:
$$
\CD
B_m(y) @>{\pi}>> B(y) \\
@V{f'}VV @VV{f}V \\
B(\tau) @>>> B(t), \\
\endCD
$$
where $B(\tau ) = \{ \tau \in \bold C \mid \mid \tau \mid < \epsilon \}$
and $\tau^m = t$.  Then $(B_m(y), \pi^{-1}(0), \pi^{-1}(Z(f)))$
is a local space and $H^i(B_m(y)-\pi^{-1}(Z(f)), \bold Q)$ has a
mixed Hodge structure defined in \S 2.1.  The covering transformation
group $\mu_m = \Gal (B(\tau) /B(t))$ acts on the mixed Hodge
structure on $H^i(B_m(y)-\pi^{-1}(Z(f)), \bold Q)$ by the functoriality.
Let $K_{MH}(\bold C, \mu_m)$ be the
Grothendieck group of mixed Hodge structures with a 
$\mu_m$-action.  Then $K_{MH}(\bold C, \mu_m)$ has a natural ring
structure arising from tensor product and graded by weight.  Since
$[\bold Q - \bold Q (1)]=[\bold Q] - [\bold Q (1)]$ 
is a non zero divisor on $K_{MH}(\bold C, \mu_m)$,
the natural homomorphism from $K_{MH}(\bold C, \mu_m)$
to the localization $K_{MH}(\bold C, \mu_m)_{loc}$ of
$K_{MH}(\bold C, \mu_m)$ with respect to $[\bold Q - \bold Q(1)]$ is
injective. (See [D-L].)

  The element $\Psi_{f,m}(\bold Q) = [H^*(B_m(y) -\pi^{-1}(Z(f)), \bold Q)]$
in $K_{MH}(\bold C, \mu_m)$ is defined by 
$\sum_{i=0}^{2\dim B(x)}(-1)^i[H^i(B_m(y) - \pi^{-1}(Z(f)), \bold Q)]$.
We define 
$$
\align
&
\phi_{f,m}(\bold Q) =
\frac{1}{[\bold Q -\bold Q(-1)]}\Psi_{f, m}(\bold Q)
\in K_{MH} ( \bold C, \mu_m)_{loc}, \text{ and } \tag{3.2.1}\\
& \Phi_{f,m}(\bold Q) =
\phi_{f,m}(\bold Q) - [\bold Q]
\in K_{MH} ( \bold C, \mu_m)_{loc}.	 \\
\endalign
$$
Note that $\phi_{f,m}(\bold Q)$ corresponds to the cohomology
of the Milnor fiber if the exponent of monodromy divides $m$.  Using
the weight spectral sequence, one can show $\phi_{f,m}$ is actually
an element in $K_{MH}(\bold C, \mu_m)$, but we do not use this fact
in the rest of this paper.

Let  $y_i = (y_{i1}, \dots , y_{im_i})$ be sets of coordinates
as in the end of \S 2.2.
We  define 
$B(y_i) = \{(y_{i1}, \dots ,y_{i,m_i}) \in \bold C^{m_i}\mid
\mid y_{ij} \mid < \epsilon \qquad (1 \leq j \leq m_i)\}$.
Let $g_i = g_i(y_i)$ be a holomorphic function on $B(y_i)$ with
$g_i(0) =0$.  Then $g_i$ defines a homomorphism $B(y_i) \to B(x_i)$
by $x_i = g_i(y_i)$.  We assume that $g_i$ is smooth on 
$B(y_i) - g_i^{-1}(0)$.  The fiber product $g$ of $g_i$ 
from $B(y) = \prod_{i=1}^nB(y_i)$ to $B(x) = \prod_{i=1}^nB(x_i)$ 
is defined as in \S 2.2.
The sequence of positive numbers $(d_1, \dots , d_n)$ is denoted by 
$\bold d$. We define $B_{\bold d}(x)$ and $\tilde B_{\bold d}(x)$
by the following cartesian squares.
$$
\CD
B_\bold d(y) @>>> B (y) 
@.\qquad @. 
\tilde B_\bold d(y) @>>> B_\bold d(y) \\
@VVV @VVV 
@. 
@VVV @VVV \\
B(\xi) @>>> B(x), 
@. @. 
B(\tau) @>>> B(t) \\
\endCD
$$

We introduce the equivariant version of the definition (3.2.1). 
The natural map $\tilde B_\bold d(y)$
$\to B(\tau)$ is denoted by
$\tilde f(y)$.
We define $\tilde\phi_{f\circ g,m}(\bold Q)$
by
$$
\align
\tilde\phi_{f\circ g,m}(\bold Q)= &
\frac{1}{[\bold Q -\bold Q(-1)]}
[H^*(\tilde B_{\bold d}(y) - (\tilde f(y))^{-1}(0), \bold Q)] - [\bold Q] \\
&\in K_{MH} ( \bold C, G\times\mu_m)_{loc}. \\
\endalign
$$
Then it is easy to see that $\Phi_{f\circ g,m}(\bold Q) =
\tilde\Phi_{f\circ g,m}(\bold Q)^G$.
For a subset $I$ of $[1, n]$, we define $B_I(\xi)$
by $\prod_{i \notin I}B(\xi_i)$.  The fiber product of
morphism $B_I(\xi) \to B(\xi) \to B(t)$ and $B(\tau) \to B(t)$
is denoted by $B_{I,m}(\xi)$. 
Put $G_I= \prod_{i\notin I}\mu_{d_i}$.
Then $G_I\times \mu_m$ acts on $B_{I,m}(\xi)$.
The natural morphism from $B_{I,m}(\xi) \to B(\tau)$ is denoted
by $\tilde f_I(\xi)$.  We define $\tilde\phi_{f_I,m}(\bold Q)$ by
$$
\align
\tilde \phi_{f_I,m}(\bold Q) = & 
\frac{1}{[\bold Q -\bold Q(-1)]}
[H^*(B_{I,m}(\xi) - (\tilde f_I(\xi))^{-1}(0), \bold Q)]
- [\bold Q] \\
& \in K_{MH} ( \bold C, G_I\times\mu_m)_{loc}. \\
\endalign
$$
Let $\tilde B^0_{I, \bold d}=\tilde B_{I, \bold d}
-\cup_{J\supsetneq I}\tilde B_{J, \bold d}$,
$\tilde B_{\bold d}^0 = \tilde B_{\emptyset ,\bold d}^0$  and the $k_{I}$
be the natural inclusion 
$$
k_I : \tilde B^0_{I, \bold d}-(\tilde f(y))^{-1}(0) \to
\tilde B_{I, \bold d}-(\tilde f(y))^{-1}(0).
$$
Then we use the following equality.
$$
[H^*
(\tilde B_{\bold d}(y) - (\tilde f(y))^{-1}(0), \bold Q)]
=
\sum_{I \subset [1,n]}
[H^*
(\tilde B_{I,\bold d}(y) - (\tilde f(y))^{-1}(0), k_{I!}\bold Q)]
\tag{3.2.2}
$$

Now we introduce several notations related to the resolution
of singularities.
Let $b(y_i): \hat B(y_i) \to B(y_i)$ be a resolution of singularity for
$g_i^{-1}(0) \subset B(y_i)$, i.e, $b(y_i)$ is a bimeromorphic projective
morphism and $b(y_i)^{-1}(g^{-1}(0))_{red}$ is a simple normal crossing
divisor.  Let $d_i$ be the least common multiple of the multiplicities of 
the components in $b(y_i)^{-1}(g^{-1}(0))$.  Let $B(\xi_i) \to B(x_i)$
be the morphism defined by $\xi_i^{d_i} = x_i$.  The normalization
of $\hat B(y_i) \times_{B(x_i)} B(\xi_i)$ is denoted by $\tilde B(y_i)$.
The second projection $\tilde B(y_i) \to B(\xi_i)$ is denoted by 
$\tilde g_i$.  
Then the morphism $\tilde g_i$ is $\mu_{d_i}$-equivariant and
$\tilde B(y_i)$ is a quasi-smooth variety with reduced
quasi-normal crossing divisor $\tilde g_i(0)$.  
$$
\CD
\tilde B(y_i) @>>> \hat B(y_i) \\
@V{\tilde g_i}VV @VV{\hat g_i}V \\
B(\xi_i) @>>> B(x_i) \\
\endCD
$$
Let $V = \tilde B(y)=
\prod_{i=1}^n\tilde B(y_i)$ be the product of $\tilde B(y_i)$ and
the product $V \to B(\xi) =\prod_{1=1}^nB(\xi_i)$ of $g_i$ is denoted by
$h$.  The base change of $h$ by the morphism $\hat B(\xi) \to B(\xi)$,
$\tilde B(\xi) \to B(\xi)$ and $\bar B(\xi )\to B(\xi)$ are denoted by
$\hat h : \hat V \to \hat B(\xi)$,
$\tilde h : \tilde V \to \tilde B(\xi)$ and
$\bar h : \bar V \to \bar B(\xi)$.  The natural morphism 
$B(\xi ) \to B(t)$ is denoted by $f(\xi)$.
Then we have the following diagram.
$$
\CD
@. \hat V @>>> V = \prod_{i=1}^n \tilde B(y_i)\\
@. @V{\hat h}VV @VV{h}V \\
@. \hat B(\xi) @>{b(\xi)}>> B(\xi) @>{f(\xi)}>> B(t) \\
@. @AAA @AAA \\
B(\tau) @<<< \tilde B(\xi) @>>> \bar B(\xi) \\
@. @A{\tilde h}AA @AA{\bar h}A \\
@. \tilde V @>>> \bar V \\
\endCD
$$
The morphism $h: V \to B(\xi)$
is $G$-equivariant and $G\times \mu_m$ acts on $\tilde V$ over
$B(\tau)$.
Let $B^0(\xi) = \prod_{i=1}^n B^0(\xi_i)$, $D_B = B(\xi) - B^0(\xi)$
and $D_F = (f(\xi))^{-1}(0)$.  The proper transform of $D_B$ and
$D_F$ in $\hat B(\xi)$ is denoted by $\hat D_B$ and $\hat D_F$ and
the exceptional divisor for $b(\xi)$ is denoted by $\hat D_E$.
The pullbacks of $\hat D_B, \hat D_F$ and $\hat D_E$ 
in $\tilde B(\xi)$ are denoted by
$\tilde D_B$, $\tilde D_F$ and $\tilde D_E$, respectively.  
The group $G \times \mu_m$ acts on $\tilde D_E$.
We define 
divisors on $\tilde V$ as follows:
$$
 \tilde h^{-1}(\tilde D_B) = \tilde V_B, 
\tilde h^{-1}(\tilde D_F) = \tilde V_F, 
\tilde h^{-1}(\tilde D_E) = \tilde V_E
$$

Let $\tilde j_V: \tilde V-(\tilde V_F \cup \tilde V_E)
\to \tilde V$ be the natural immersion.  Then the mixed Hodge structure
of $H^*(\tilde V -(\tilde V_F \cup\tilde V_E), \bold Q)$ is given
by $\bold H^*(\tilde V, \bold R \tilde j_{V*}\bold Q \mid_{\tilde V_E})$.
The inverse image of $(0, \dots ,0)$ under the morphism
$V \to \prod_{i=1}^nB(y_i)$ is denoted by $V_c$.  
The action of $G$ on $V$ induces that of $V_c$.
The inverse image of $V_c$
under the morphism $ \tilde V \to V$ is denoted by $\tilde V_c$.
Then we have $\tilde V_c = V_c \times \tilde D_E$
and the action of $G \times \mu_m$ is equal to the diagonal action.
The restriction of $\tilde h$ to $\tilde V_c$ is denoted by 
$\tilde h_c : \tilde V_c \to \tilde D_E$.  It is easy to see that $V_c \to V$
and $\tilde V_c \to \tilde V$ are retractions.  
The group $G \times \mu_m$ acts on the mixed Hodge structure
$\bold H^i(\tilde V, \bold R \tilde j_{V*}\bold Q)
\simeq \bold H^i(\tilde V_c, \bold R \tilde j_{V*}\bold Q
\mid_{\tilde V_E})$.  Therefore we have
$$
[\bold H^*(\tilde V, \bold R \tilde j_{V*} \bold Q)]
=[H^* (\tilde V -(\tilde V_E\cup \tilde V_F), \bold Q)]
\in K_{MH}(\bold C, G \times \mu_m).
$$
Let $\tilde D_E =
\coprod_{\sigma \in F'}\tilde Z_{\sigma}^0$ be the stratification
indexed by $F'$ and $\tilde Z_{\sigma}$ be the closure of 
$\tilde Z_{\sigma}^0$ in $\tilde D_E$.  
Note that this stratification is stable under the action of 
$G \times \mu_m$.
Let $\tilde U_{\sigma}$ be a
tubular neighborhood of $\tilde Z_{\sigma}$, 
$j_{\sigma}: \tilde U_{\sigma}-(\tilde D_E\cup \tilde D_F) 
\to \tilde U_{\sigma}$
and $\tilde j_{\sigma}:
\tilde h^{-1}(\tilde U_{\sigma}-(\tilde D_E \cup \tilde D_F)) \to
\tilde h^{-1}(\tilde U_{\sigma})$ be the natural inclusion.
$$
\CD
\tilde h^{-1}(\tilde U_{\sigma}-(\tilde D_E \cup \tilde D_F)) 
@>{\tilde j_{\sigma}}>>
\tilde h^{-1}(\tilde U_{\sigma}) 
@<{\supset}<<  \tilde h_c(\tilde Z_\sigma^0) \\
@VVV @VVV @VVV \\
\tilde U_{\sigma}-(\tilde D_E \cup \tilde D_F)
@>{j_{\sigma}}>>
\tilde U_{\sigma} @<<{\supset}< \tilde Z_\sigma^0 \\
\endCD
$$
Let $G_{\sigma}$ be the stabilizer of $\tilde Z^0_{\sigma}$
in $G\times \mu_m$.  Then 
$\bold H^i_c(\tilde h^{-1}_c(\tilde Z_{\sigma}^0), 
\bold R \tilde j_{V*} \bold Q)$
is a $G_{\sigma}$-module and we have
$$
\align
[\bold H^*(\tilde V_c, \bold R \tilde j_{V*} \bold Q)]
& =
\sum_{[\sigma] \in F'/(G\times \mu_m)}
\Ind_{G_\sigma}^G[
\bold H^*_c(\tilde h^{-1}_c(\tilde Z_{\sigma}^0), 
\bold R \tilde j_{V*} \bold Q)] \\
& =
\sum_{[\sigma] \in F'/(G\times \mu_m)}
\Ind_{G_\sigma}^G[
\bold H^*_c(\tilde h^{-1}_c(\tilde Z_{\sigma}^0), 
\bold R \tilde j_{\sigma *} \bold Q)]. \\
\endalign
$$
Put $I(\sigma )=\{i \mid \bold e_i \in \sigma \}$ and
$\Psi_{V, m,\sigma}(\bold Q)=
[\bold H^*_c(\tilde h^{-1}_c(\tilde Z_{\sigma}^0), 
\bold R \tilde j_{\sigma *} \bold Q)]$.
Then we have
$$
[H^*(\tilde B_\bold d(y) - (\tilde f (y))^{-1}(0), 
k_{\emptyset !}k_{\emptyset}^*\bold R\tilde j_\sigma \bold Q)]
=\sum_{\{\sigma \in F', I(\sigma) = \emptyset\}/(G \times \mu_m)}
\Ind_{G_\sigma}^G(
\Psi_{V, m,\sigma}(\bold Q))
$$
in $K_{MH}(\bold C, G \times \mu_m)$.
In the next section, we compute $R^i\tilde j_{\sigma *}\bold Q$
using toric geometry.

\heading
\S 3.3 The structure of $R^i \tilde j_{\sigma *}\bold Q$.
\endheading

First we describe a stratification of $\tilde D_E$ in terms of toric
geometry.  The space $B^0(\xi, \tau)=B^0(\xi) \times B^0(\tau)$ can be 
considered as an open set of $\bar B(\xi)$.  The complement
$\bar D = \bar B(\xi) - B^0(\xi) \times B^0(\tau)$ is a 
quasi normal crossing divisor whose irreducible components are indexed by
the set of 1-dimensional cones $r$ in $\tilde F$.
(See (1.2.1) for the definition of $\tilde F$.)  
Let $\{ r_0, r_1, \dots ,r_s, r_{s+1}, \dots , r_{s+n}\}$ 
be the set of 1-dimensional cones of 
$\tilde F$ such that $r_0,r_{s+1}, \dots , r_{s+n}$ correspond to
the proper transforms of $\{\tau = 0\}$, $\{ \xi_1 = 0\}, \dots ,
\{ \xi_n =0\}$, respectively.  The corresponding divisors are written as
$\bar D_0, \bar D_1, \dots ,\bar D_s, \bar D_{s+1}, \dots , \bar D_{s+n}$.
The intersection
$\bar D_{i_1} \cap \dots \cap \bar D_{i_k}$ is non-empty and as a consequence
defines a stratum of $\bar D$ if and only $r_{i_1}, \dots ,r_{i_k}$ generate
a cone $\sigma$ of the fan $\tilde F$.  We denote $\bar Z_{\sigma} =
\bar D_{i_1} \cap \dots \cap \bar D_{i_k}$.  
Let $\bar D_F =\bar D_0$, $\bar D_B =\cup_{i=1}^n\bar D_{s+i}$
and $\bar D_E = \cup_{i=1}^s\bar D_i$.
Since 
$\bar D_F \cap \tilde B(\xi ) = \tilde D_F$, 
$\bar D_B \cap \tilde B(\xi) =\tilde D_B$ and
$\bar D_E \cap \tilde B(\xi) = \tilde D_E$,
$\tilde Z_{\sigma} = \bar Z_{\sigma} \cap \tilde B(\xi)$, where $\sigma$
is a simplicial cone of $\tilde F$, defines a stratification of 
$\tilde D_F \cup \tilde D_E \cup \tilde D_B$.  It is easy to see that
$\tilde Z_{\sigma}$ is not empty if and only if $\dim \bar Z_{\sigma} \geq 1$
i.e. $\dim \sigma < n+1$.  We put $\bar Z_{\sigma}^0 =
\bar Z_{\sigma} - \cup_{\tau > \sigma}\bar Z_{\tau}$.  Then 
$\bar Z_{\sigma}^0$ is isomorphic to a torus and $\tilde Z_{\sigma}^0$
is a non-singular hypersurface in $\bar Z_{\sigma}^0$.  As a consequence,
we have a stratification
$$
\tilde D =
\tilde D_B \cup \tilde D_E \cup \tilde D_F =
\cup_{\sigma:\text{ simplicial }, \dim \sigma < n+1}\tilde Z_{\sigma}^0.
$$
Under this stratification, $\tilde D_E$ corresponds to $\sigma$'s
whose generator contains at least one of $r_1, \dots , r_s$.  Let
$A_{\sigma} = \bold C[\hat \sigma \cap L(\xi, \tau)]$ and 
$Spec(A_{\sigma}) \to X_{\tilde F} $ be the natural map.  Let $\bar U_{\sigma}$
be the pull back of $Spec(A_\sigma)$ in $\bar B(\xi)$.  Let $L_{\hat\sigma}$
be the maximal linear subspace contained in $\hat \sigma$ and
$\bar A_{\sigma} = \bold C[L_{\hat\sigma} \cap L(\xi,  \tau)]$.
The kernel of the natural map $A_{\sigma} \to \bar A_{\sigma}$ is denoted by
$I_{\sigma}$.  Then the corresponding closed subvariety is identified
with $\bar Z_{\sigma}^0$.
$$
\CD
\bar Z_{\sigma}^0 @>{\subset}>> \bar U_{\sigma} @>>> \bar B(\xi) \\
@. @V{\cap}VV @VV{\cap}V \\
@. Spec(A_{\sigma}) @>>> X_{\tilde F} \\
\endCD
$$
By taking the intersections with $\tilde B(\xi)$, we have the following
diagram.
$$
\CD
\tilde h_c^{-1}(Z_{\sigma}^0) @>{\subset}>> \tilde h^{-1}(\tilde U_{\sigma})
@<{\tilde j_{\sigma}}<< \tilde h^{-1}
(\tilde U_{\sigma}-(\tilde D_E\cup \tilde D_F)) \\
@VVV @VVV @VVV \\
\tilde Z_{\sigma}^0 @>>{\subset}> \tilde U_{\sigma}
@<<{\supset}< \tilde U_\sigma -(\tilde D_E \cup \tilde D_F),
\endCD
$$
where $\tilde U_{\sigma} = \bar U_{\sigma}\cap\tilde B(\xi)$.
We use these varieties to compute 
\linebreak
$\bold H^*(\tilde h_c^{-1}(Z_{\sigma}^0), \bold R \tilde j_{\sigma}
\bold Q\mid_{h_c^{-1}(\tilde Z_{\sigma}^0)})$.
Let $\bar V_E = \bar h^{-1}(\bar D_E)$, $\bar V_F = \bar h^{-1}(\bar D_F)$
and $\bar V_B = \bar h^{-1}(\bar D_B)$.
The natural morphisms $\tilde V \to V$ and $\tilde V \to \tilde B(\xi)$ 
define a closed immersion $\tilde V \to V \times \tilde B(\xi)$.
This morphism defines a morphism of local spaces:
$$
\align
(\tilde V, \tilde V_c , \tilde V_E\cup \tilde V_F \cup \tilde V_B)
& \overset{\tilde\alpha}\to\to 
(V \times \tilde B(\xi), V_c \times \tilde D_E, 
((V-V^0) \times \tilde B(\xi)) \cup (V \times \tilde D)), \\
(\bar V, \bar V_c, \bar V_E \cup \bar V_F \cup \bar V_B)
& \overset{\bar\alpha}\to\to
(V \times \bar B(\xi) , V_c \times \bar D_E, V \times \bar B(\xi)-
V^0\times B^0(\xi, \tau)),\\
\endalign
$$
where $\tilde D = \tilde D_E \cup \tilde D_F \cup \tilde D_B$ and
$V^0 = h^{-1}(B^0(\xi))$.  Note that $V_c \times \tilde D_E \simeq
\tilde V_c$ and $V_c \times \bar D_E = \bar V_c$.
The inclusion of divisors $\tilde D_E \cup \tilde D_F \to
\tilde D_E \cup \tilde D_F \cup \tilde D_B$ defines the following 
open immersion of local spaces:
$$
\align
(\tilde V, \tilde V_c , \tilde V_E\cup \tilde V_F \cup \tilde V_B)
&\overset{\tilde \beta}\to\to 
(\tilde V, \tilde V_c , \tilde V_E\cup \tilde V_F ), \\
(\bar V, \bar V_c, \bar V_E \cup \bar V_F \cup \bar V_B)
& \overset{\bar\beta}\to\to
(\bar V, \bar V_c, \bar V_E \cup \bar V_F ). \\
\endalign
$$ 
We investigate these inclusions on the open set 
$\tilde h^{-1}(\tilde U_{\sigma}) \subset \tilde V$ and
$V \times \tilde U_{\sigma}$.

Let $y$ be a point in $\tilde B(y)$ and $z_{11}, \dots , z_{1m_1}, \dots ,
z_{n1}, \dots , z_{nm_n}$ be a local coordinate near $y$ such that
$g_1, \dots , g_n$ can be written as
$$
g_1 = z_{11}^{m_{11}}\cdots z_{1l_1}^{m_{1l_1}},
\dots, g_n = z_{n1}^{m_{n1}}\cdots z_{nl_n}^{m_{nl_n}}.
$$
Then the local equation of $V$ at a lifting $\tilde y$ of $y$ 
can be written as
$$
\xi_1^{d_1'}  = z_{11}^{m_{11}'}\cdots z_{1l_1}^{m_{1l_1}'}, \dots , 
\xi_n^{d_n'} = z_{n1}^{m_{n1}'}\cdots z_{nl_n}^{m_{nl_n}'}, 
$$
where 
$$
d_i' = \frac{1}{\gcd (d_i, m_{i,j})_j}d_i,
m_{ij}' = \frac{1}{\gcd (d_i, m_{i,j})_j}m_{ij}.
$$

  Let $r_1, \dots , r_a$ be a generator of $\sigma \cap (L(\xi, \tau) \otimes
\bold Q)$, such that 
$\sum_{i=1}^a\bold N r_i \subset (\sigma \cap L(\xi, \tau))$.  
Suppose that $r_1, \dots , r_k$ correspond to the components of $\bar D_B$
and $r_{k+1}, \dots ,$
\linebreak
$r_a$ correspond to the components of 
$\bar D_E \cup \bar D_F$.
We take 
$r_{a+1}, \dots , r_{n+1}$ such that $\sum_{i=1}^{n+1}\bold Z r_i$
is finite index in $L(\xi, \tau)^*$.  
Then by changing the numbering of $\xi_i$, 
$r_1, \dots , r_a$ can be written as
$$
\align
r_1 &= (r_{11}, 0, \dots ,0), \dots, 
r_k = (0, \dots ,r_{kk} , 0, \dots , 0), \\ 
r_{k+1} & = (r_{k+1,1}, \dots , r_{k+1,n}, r_{k+1, n+1}),\dots , 
r_a  = (r_{a1}, \dots , r_{an}, r_{a, n+1}) \\
\endalign
$$
Then $\bold C[\hat \sigma \cap L(\xi, \tau )]$ is a subring of
$\bold C [u_1, \dots, u_a, u_{a+1}^{\pm}, \dots , u_{n+1}^{\pm}]$,
where $\{ u_i\}_i$ is the multiplicative 
expression of the dual base $\{ r_i^*\}_i$ of $\{r_i\}_i$.
The subring $\bold C[\hat\sigma \cap L(\xi, \tau)]$ is 
characterized by an invariant subring of 
$\bold C[u_1, \dots, u_a, u_{a+1}^{\pm}, \dots , u_{n+1}^{\pm}]$
under an action of an abelian group $H$.
The morphism 
$Spec( \bold C[u_1, \dots, u_a, u_{a+1}^{\pm}, \dots , u_{n+1}^{\pm}] \to
Spec(\bold C [\xi_1, \dots ,\xi_n])$ is given by
$$
\align 
\xi_1 & = u_1^{r_{11}} \cdots u_a^{r_{a1}} \cdot (\text{ unit }), \dots , 
\xi_n  = u_1^{r_{1n}} \cdots u_a^{r_{an}} \cdot (\text{ unit }), \\
\tau & = u_1^{r_{1,n+1}} \cdots u_a^{r_{a,n+1}} \cdot (\text{ unit }) \\
\endalign
$$

First we investigate the inclusions $\bar \alpha$ and $\bar \beta$.
The local equation for $\bar V$ is given by
$$
\align
[ u_1^{r_{11}}\prod_{p=k+1}^a u_p^{r_{p1}}]^{d_1'}\cdot (\text{ unit })
& =z_{11}^{m_{11}'}\cdots z_{1l_1}^{m_{1l_1}'} \\
&\cdots \\
[ u_k^{r_{kk}}\prod_{p=k+1}^a u_p^{r_{pk}}]^{d_k'}\cdot (\text{ unit })
& =z_{k1}^{m_{k1}'}\cdots z_{kl_k}^{m_{kl_k}'} \\
[\prod_{p=k+1}^a u_p^{r_{p,k+1}}]^{d_{k+1}'}\cdot (\text{ unit })
& =z_{k+1,1}^{m_{k+1,1}'}\cdots z_{k+1,l_{k+1}}^{m_{k+1,l_{k+1}}'} \\
&\cdots \\
[\prod_{p=k+1}^a u_p^{r_{pn}}]^{d_n'}\cdot (\text{ unit })
& =z_{n1}^{m_{n1}'}\cdots z_{nl_n}^{m_{nl_n}'}, \\
\endalign
$$
using the coordinate of the covering $Spec(\bold C[u_1, \dots , u_{a+1},
u_{a+1}^{\pm}, \dots , u_{n+1}^{\pm}]) \times \tilde B(y)$ of
$Spec(\bold C[\hat \sigma \cap L(\xi, \tau)])\times \tilde B(y)$.
Moreover, the open set corresponding to $\bar V -(\bar V_E\cup \bar V_F \cup
\bar V_B)$ and $\bar V- (\bar V_E \cup \bar V_F)$ is given by
$$
u_1 \in B(u_1)^0 , \dots , u_a \in B(u_a)^0, \text{ and }
$$
$$
u_{k+1} \in B(u_{k+1})^0 , \dots , u_a \in B(u_a)^0,
$$
respectively.  Let 
$\bar j_{EFB} :\bar V -(\bar V_E\cup \bar V_F \cup \bar V_B) \to \bar V$,
$\bar j_{EF} :\bar V -(\bar V_E\cup \bar V_F ) \to \bar V$ and
$\bar j_M: V^0 \times B^0(\xi, \tau) \to V \times \bar B(\xi)$
be the natural inclusions.  Then
$$
\align
R^i(\bar j_{EFB})_*\bold Q & =
\wedge^i R^1(\bar j_{EFB})_*\bold Q, \tag{3.3.1} \\
R^i(\bar j_{EF})_*\bold Q & =
\wedge^i R^1(\bar j_{EF})_*\bold Q,  \\
R^i(\bar j_{M})_*\bold Q & =
\wedge^i R^1(\bar j_{M})_*\bold Q,  \\
\endalign
$$
for $i \geq 1$.
Using the local expression of $\bar V- (\bar V_E \cup \bar V_F \cup
\bar V_B)$, 
$\bar V- (\bar V_E \cup \bar V_F)$ and
$V^0\times B(\xi, \tau)^0 \to V \times \bar B(\xi )$, the
stalks if the natural homomorphism
$$
R^1(\bar j_M)_*\bold Q
\overset{\bar \alpha^*}\to\to
R^1(\bar j_{EFB})_*\bold Q
\overset{\bar \beta^*}\to\longleftarrow
R^1(\bar j_{EF})_*\bold Q
$$
at $\bar v$ are computed as follows.
$$
\align
(R^1(\bar j_M)_*\bold Q)_{\bar v} &=
[\bold Q^a \oplus( \oplus_{i=1}^n \bold Q^{l_i} \oplus \bold Q)]
\otimes \bold Q(-1) \\
(R^1(\bar j_{EFB})_*\bold Q)_{\bar v} &=
[\Coker (\bold Q^{n+1} \overset{A_{EFB}}\to\to 
(\bold Q^a \oplus( \oplus_{i=1}^n \bold Q^{l_i} \oplus \bold Q)))]
\otimes \bold Q(-1) \\
(R^1(\bar j_{EF})_*\bold Q)_{\bar v} &=
[\Coker (\bold Q^{n-k+1} \overset{A_{EF}}\to\to 
(\bold Q^{a-k} \oplus( \oplus_{i=k+1}^n \bold Q^{l_i} \oplus \bold Q)))]
\otimes \bold Q(-1), \\
\endalign
$$
where
$$
\align
&A_{EFB}(e_i)  = 
d_i'(r_{1i}, \dots, r_{ai}) \oplus (0,\dots ,0,(m_{i1}',\dots,
m_{i,l_i}'), 0, \dots , 0), 
 (1\leq i \leq n) \\
&A_{EFB}(e_{n+1})  =
d_{n+1}'(r_{1,n+1}, \dots, r_{a,n+1}) \oplus (0,\dots ,0,1),
\\
\endalign
$$
and
$$
\align
& A_{EF}(e_i)  = 
d_i'(r_{k+1i}, \dots, r_{ai}) \oplus (0,\dots ,0,(m_{i1}',\dots,
m_{i,l_i}'), 0, \dots , 0), 
 (k+1 \leq i \leq n) \\
& A_{EF}(e_{n+1})  = 
d_{n+1}'(r_{k+1,n+1}, \dots, r_{a,n+1}) \oplus (0,\dots ,0,1).
\\
\endalign
$$
The morphism $\bar \alpha^*$ is identified with the natural projection.
Since the diagram
$$
\CD
\bold Q^{n+1}@>{A_{EFB}}>>
\bold Q^a \oplus( \oplus_{i=1}^n \bold Q^{l_i} \oplus \bold Q) \\
@AAA @AAA \\
\bold Q^{n-k+1}@>>{A_{EF}}>
\bold Q^{a-k} \oplus( \oplus_{i=k+1}^n \bold Q^{l_i} \oplus \bold Q) \\
\endCD
$$
is commutative, $\bar \beta^*$ is identified with the homomorphism
induced by the natural inclusions.
By (3.3.1) and the expression of $\bar \alpha^*$ and $\bar \beta^*$,
we have the following proposition.
\proclaim{Proposition 3.3.1}
\roster
\item
The morphism $\bar\alpha^*$ is injective an $\bar\beta^*$
is surjective.
\item
The morphism $R^i(\bar j_{M})_*\bold Q \to R^i(\bar j_{EFB})_*\bold Q$
and
$R^i(\bar j_{EF})_*\bold Q \to R^i(\bar j_{EFB})_*\bold Q$
are identified with $\wedge^i \bar\alpha^*$ and
$\wedge^i \bar\beta^*$.  Moreover they are surjective and injective 
respectively.
\endroster
\endproclaim

Since $\tilde B(\xi)$ meets $\bar D_E \cup 
\bar D_F \cup \bar D_B$ transversally, we get the following proposition.
\proclaim{Proposition 3.3.2}
For $\bar v \in \tilde V$, we have
$$
\align
(R^i(\tilde j_{EFB})_*\bold Q)_{\bar v} & =
(R^i(\bar j_{EFB})_*\bold Q)_{\bar v}, 
(R^i(\tilde j_{EF})_*\bold Q)_{\bar v}  =
(R^i(\bar j_{EF})_*\bold Q)_{\bar v}, \\
(R^i(\tilde j_{M})_*\bold Q)_{\bar v} & =
(R^i(\bar j_{M})_*\bold Q)_{\bar v}. \\
\endalign
$$
As a consequence, Proposition 3.3.1 holds by replacing
$\bar j_{EFB}$, $\bar j_{EF}$ and $\bar j_M$ by
$\tilde j_{EFB}$, $\tilde j_{EF}$ and $\tilde j_M$.
\endproclaim

\heading
\S 3.4 Supplement for mixed Hodge structures
\endheading

Let $X$ be a topological space and $K=((K_{\bold Q}, W),(K_{\bold C}, W,F))$
and $K'=$\linebreak
$((K'_{\bold Q}, W),(K'_{\bold C}, W,F))$ be cohomological
mixed Hodge complexes (=CMHC's for short).  A morphism from $K$ to $K'$
is defined by the pair of morphisms of filtered and bifiltered
complexes $\varphi_{\bold Q} : K_{\bold Q}\to K'_{\bold Q}$
and $\varphi_{\bold C}: K_{\bold C} \to K'_{\bold C}$ where 
$(\varphi_{\bold Q}\otimes \bold C,W)$ and $(\varphi_{\bold C},W)$ 
are filtered congruent to each other.
After [Dur], we define the cone $\Cone (\varphi )$ of $\varphi$ by
$$
\align
\Cone (\varphi)_A^p &= K_A^p\oplus 
(K')_A^{p-1}, \\
d: \Cone (\varphi )_A^p \to \Cone (\varphi )_A^{p+1} &;
(x, y) \mapsto (dx, \varphi (x) + dy), \\
W_k\Cone (\varphi )_A^p &= W_kK_A^p \oplus W_{k+1}(K')_A^{p-1}, \text
{ for $A = \bold Q, \bold C$ 
and }\\
F^q\Cone (\varphi )_{\bold C}^p 
& = F^qK_{\bold C}^p \oplus F^q(K')_{\bold C}^{p-1}.
\endalign
$$
According to [Dur], $\Cone (\varphi )$ give rise to be a CMHC.
\demo{Definition 3.4.1}
A morphism $\varphi : K \to K'$ of CMHC's is called weakly equivalent
if the underlying morphism $K_{\bold Q} \to K'_{\bold Q}$ is a
quasi-isomorphism.
\enddemo
\demo{Remark 3.4.2}
By weak equivalence, we do not impose that $\varphi_{\bold Q}$
or $\varphi_{\bold C}$ are filtered quasi-isomorphism.
Therefore, in general, the spectral sequences 
$$
\align
& E_1^{p,q}(K,W)=\bold H^{p+q}(X, Gr^W_{-p}K) \Longrightarrow
E_\infty^{p+q}(K,W) = \bold H^{p+q}(X, K) \\
& E_1^{p,q}(K',W)=\bold H^{p+q}(X, Gr^W_{-p}K') \Longrightarrow
E_\infty^{p+q}(K',W) = \bold H^{p+q}(X, K') \\
\endalign
$$
does not coincide.  But they degenerate at $E_2$-terms and the filtration
induced by this spectral seequence is equal to weight filtration,
they coincides at $E_2$-terms.
\enddemo
\demo{Example 3.4.3}
Let $X$ be a smooth algebraic variety and $X_1$ and $X_2$ be a smooth 
comactification with normal crossing boudaries $D_1$ and $D_2$.  
Assume that there
exists a morphism $f:X_1 \to X_2$ which induces an identity on $X$.
Let $j_i : X \to X_i$ be a natural inclusion for $i =1,2$.
As in [Del], 
$K_i =((\bold Rj_{i *}\bold Q, \sigma ),
(\Omega_{X_i}^{\bullet}(\log D_i), 
W, F))$
is a CMHC on $X_i$ for $i =1,2$. Therefore
$$
\bold Rf_* K_1 =((\bold Rf_*\bold Rj_{1 *}\bold Q, \bold Rf_*(\sigma )),
(\bold Rf_*\Omega_{X_1}^{\bullet}(\log X), 
\bold Rf_*(W), \bold Rf_*(F)))
$$
is a CMHC on $X_2$ and there is a natural morphism $K_2 \to \bold Rf_*K_1$.
In general, they are not filtered quasi-isomorphic, but they are
weakly equivalent.
\enddemo
The following lemma will be useful later.
\proclaim{Lemma 3.4.4}
\roster
\item
If $f:K \to K'$ is weakly equivalent, then the homomorphism
$H^i(f) :\bold H^i(X, K) \to \bold H^i(X, K')$ is an isomorphism of 
mixed Hodge strucutres.
\item
Let $f:K \to L$ and $b:L' \to L$ be morphisms of CMHC's on $X$.
Suppose that $b$ is weak equivalent.  Then there exists a cohomological
mixed Hodge complex $K'$ and morphisms $a: K' \to K$ and 
$f':K' \to L'$ of CMHC's on $X$ such that $a$ is a weak equivalence and
$f\circ a = b \circ f'$.
\endroster
\endproclaim
\demo{Proof}

(1) The homomorphism
$H^i(f)$ preserves the weight and Hodge filtraions and compatible 
with the complex conjugate and bijective.  Since the category of
mixed Hodge structures is an abelian category, it is an isormophims.

(2) Let $K'= \Cone ((f,b):K \oplus L' \to L)$ and 
$a :K' \to K$ and $f':K' \to L'$ be the natural morphisms.  Then $K'$
is equipped with CMHC and $a$ and $f$ are morphism of CMHC's which
satisfy the conditions of the lemma.
\enddemo
We define direct sum decomposition of CMHC's.
\demo{Definition 3.4.5}
\roster
\item
Let $K$, $L$ and $M$ be CMHC's and $p_1: K \to L$ and 
$p_2 :K \to M$ be morphisms of CMHC's.  $p_1 \oplus p_2 : K \to L \oplus
M$ is called a direct sum decomposition if it is weakly equivalent.
\item
Two direct sum decompositions $f_1: K\to L_1 \oplus M_1$ and
$f_2 : K \to L_2 \oplus M_2$ are said to be equivalent if there exists 
a third direct sum decomposition $p_3 : K \to L_3 \oplus M_3$
and weak equivalences $g_i : L_3 \to L_i$ and $h_i : M_3 \to M_i$
$(i=1,2)$ such that the following diagram commutes
$$
\CD
K @>{f_3}>> L_3 \oplus M_3 \\
@V{=}VV @VV{g_i\oplus h_i}V \\
K @>>{f_i}> L_i \oplus M_i \\
\endCD
$$
\endroster
\enddemo
\proclaim{Lemma 3.4.6}
Let $f:K \to L \oplus M$ be a direct sum decomposition and $K \to K'$
be a weak equivalence.  Then there exist
CMHC's $L'$, $M'$ and weak equivalences
$a: L \to L'$, $b: M \to M'$ and $f': K' \to L' \oplus M'$
such that the following diagram commutes.
$$
\CD
K @>{g}>> K' \\
@V{f}VV @VV{f'}V \\
L \oplus M @>>{(a,b)}> L' \oplus M' \\
\endCD
$$
\endproclaim
\demo{Proof}
Since $f$ is weakly equivalent, the natural morphism
$f^*:\Cone (K \to L) \oplus \Cone (K\to M) \to K$
is weakly equivalent.  Therefore the composite 
$$
g \circ f^* : \Cone (K \to L) \oplus \Cone (K \to M) \to K'
$$
is weakly equivalent.  Therefore
$$
f' :K' \to \Cone (\Cone (K \to L) \to K')[1] \oplus
\Cone (\Cone (K \to M) \to K') [1]
$$
is a direct sum decomposition.  We put $L' =
\Cone (\Cone (K \to L) \to K') [1]$ and
$M' = \Cone (\Cone (K \to M) \to K') [1]$.
Since the natural morphisms $L \to \Cone (\Cone (K \to L) \to K)$
and $M \to \Cone (\Cone (K \to M) \to K)$ are weakly equivalent,
the composites
$$
\align 
a: L  &  \to \Cone (\Cone (K \to L) \to K)
\to \Cone (\Cone (K \to L) \to K')= L' \\
b: M  & \to \Cone (\Cone (K \to M) \to K)
\to \Cone (\Cone (K \to M) \to K')= M' \\
\endalign
$$
are weakly equivalent.  Thus we have the required CMHC's and weak
equivalences.
\enddemo
\demo{Definition 3.4.7 (Quasi-canonical filtration)}
\roster
\item
Let $K$ be a CMHC on $X$.  A sequence $\kappa_i:K_i \to K$ and 
$\kappa_{i,i+1}:K_i \to K_{i+1}$
of morphism of CMHC's is called quasi-canonical filtration if 
(1)
there exists $K'_i \to K$ such that the diagram is commutative
$$
\CD 
K'_i @>{a}>> K_i \\
@V{b}VV @VVV \\
\tau_iK @>>> K, \\
\endCD
$$
where $a$ and $b$ are quasi-isomorphism.  Note that we do not impose
that either $K'_i$ or $\tau_iK$ are equipped with structures of
CMHC's, and (2)
$\kappa_i = \kappa_{i+1}\circ \kappa_{i,i+1}$.
\item
The decomposition $\Cone (K_{i-1} \to K_i)[1] \to L_i \oplus M_i$
of the cone $\Cone (K_{i-1}$ \linebreak
$\to K_i)$ is called a decomposition associated
to the quasi-canonical filtration $\{ K_i \to K, K_{i-1} \to K_i\}$.
\endroster
\enddemo
\demo{Definition 3.4.8 (Associated decomposition)}
Let $K_i \to K, K_{i-1} \to K_i$ be a quasi-canonical filtration
and $\Cone (K_{i-1} \to K_i)[1] \to C^i \oplus I^i$ be a decomposition
associated to the quasi-canonical filtration.  Let $L \to K$ be a morphism
of CMHC's.  The decomposition is associated to the morphism $L \to K$ if
there exists complexes $A$, $B$ and morphisms $f:A \to B$,
$a:A \to \Cal H^i L$, $b: B \to \Cal H^i K$ and 
$b': B \to \Cone (K_{i-1} \to K_i)[1]$ such that
\roster
\item
The morphisms $a$, $b$ and $b'$ are quasi-isomorphisms.
\item
The following diagram commutes
$$
\CD
A @>{f'}>> B @>{b'}>> \Cone (K_{i-1} \to K_i)[1]\\
@V{a}VV @VV{b}V \\
\Cal H^i L @>>{\Cal H^i(f)}> \Cal H^i K. \\
\endCD
$$
\item
The composite $A \to B \to \Cone (K_{i-1}\to K_i)[1] \to I^i$ 
is a quasi-iso-
\linebreak
morphism.
\endroster
\enddemo
\demo{Remark 3.4.9} We do not impose that 
the complexes of sheaves $A$, $B$ and
$\Cal H^i L$ are equipped with the structure of CMHC's.
\enddemo
Now we can prove the following proposition.
\proclaim{Proposition 3.4.10}
Let $L$, $K$ be CMHC's, $K_i \to K, K_{i-1} \to K_i$ be
a finite quasi-canonical filtration, $f_i :\Cone (K_{i-1} \to K_i)[1]
\to C^i \oplus I^i$ be decompositions of $\Cone (K_{i-1} \to K_i)$ 
and $\varphi :L \to K$ be a morphism of CMHC's.
Suppose that the morphism $\varphi$ is associated to
the decompositions $f_i$.
Then there exists a quasi-canonical filtration $\{L_i \to L,
L_{i-1} \to L_i\}$ such that $\Cone (L_{i-1} \to L_i)$ is
weakly equivalent to $I^i$.
\endproclaim
\demo{Proof}
We construct $L_m$, $L_m \to K_m$ and $L_{m-1} \to L_m$ by descending
induction.  For a sufficiently large $m$, $K_m \to K$ is a quasi 
isomorphism and $K_i \to K_{i+1}$ is isomorphic for $i \geq m$.
By Lemma 3.4.4, we can take a CMHC $L^m$
such that the following  diagram of CMHC's
commutes
$$
\CD
L_m @>>> K_m \\
@V{i_m}VV @VVV \\
L @>>> K \\
\endCD
$$
such that $i_m$ is a weak equivalence.  
Therefore we can take $L_m \to L$ for a sufficiently large $m$.
Suppose that $L_m$, $L_m \to K_m$ and $L_m \to L_{m+1}$
$(m \geq k)$ satisfying the conditions of the proposition are given
by induction.  
Let $'I^{k} =
\Cone (\Cone (K_{k-1} \to K_k)[1] \to C^k)$.  Then the composite
$'I^k \to \Cone (K_{k-1} \to K_k)[1] \to I^k$ is a weak equivalence
and it is denoted by $b$.
The morphism $L_k \to K_k \to \Cone (K_{k-1} \to K_k)[1]
\to I^k$ is denoted by $a$.  Then we have the following diagram;
$$
\CD
@. 'I^k @>>> \Cone (K_{k-1} \to K_k)[1] \\
@. @VV{b}V \\
L^k @>{a}>> I_k, \\
\endCD
$$
where $b$ is weakly equivalent.  Then by Lemma, there exists $L'_k$,
$a'' : L'_k \to 'I^k$ and $b': L'_k \to L_k$ such that $b'$ is weakly
equivalent.  Let $L_{k-1} = \Cone (L'_k \to 'I^k)$.  Then there
are natural morphisms
$$
L_{k-1} = \Cone (L'_k \to 'I^k) \to 
\Cone (K_k \to \Cone (K_{k-1} \to K_k)[1])
\to K_{k-1}
$$
and $L_{k-1} = \Cone (L'_k \to 'I^k) \to L'_k \to L_k$.
This proceeds the steps of the induction.
\enddemo
\proclaim{Corollary 3.4.11}
Under the notations as in Proposition 3.4.10,
$$
[\bold H^*(X, L)] = \sum_{i} (-1)^i [\bold H^*(X, I^i)]
$$
in $K_{MH}(\bold C)$
\endproclaim
\demo{Proof}
Since $L$ is finite, $L_m$ is exact for $m << 0$ and $L_m \to L$
is weakly equivalent. for $m >> 0$.  Therefore, we have 
$$
\align
[\bold H^*(X, L)] &= \sum_{i}(-1)^i [\bold H^*(X , \Cone (L_{i-1} \to L_i))] \\
&
=\sum_i (-1)^i[\bold H^* (X, I^i)].
\endalign
$$
\enddemo

\heading
\S 3.5  Direct sum decomposition for normal crossing local spaces
\endheading

  Let $(N, X, Y)$ be a normal crossing local space and $\cup_{i \in I} D_i$
be the irreducible decomposition of $D = X \cup Y$.  Let
us define a subvariety $D_J$ as $\cap_{j \in J} D_j$.
Suppose that $D_J \subset X$.
The complex of sheaves $\Omega_N^{\bullet}(\log D) \otimes \Cal O_{D_J}$ 
is denoted by $\Omega^{\bullet} (J)$.  We have introduced a structure of
cohomological mixed Hodge complex on $\bold R j_* \bold Q\mid_{D_J}$,
where $j : N- (X\cup Y) \to N$ is the natural inclusion, via the
quasi-isomorphism
$$
(\bold Rj_* \bold Q\mid_{D_J})\otimes \bold C \simeq \Omega^{\bullet}(J)
$$
with the filtrations $W_k$ and $F^p$ on the right hand side.  Let 
$D_J^0 = D_J - \cup_{i \notin J}D_{J\cup i}$ and
the natural inclusion from $D^0_J$ to $D_J$ is denoted by
$j_J : D_J^0 \to D_J$.  As in Corollary 3.1.10, the structure of CMHC on
$(j_J)_!(\bold R j_*\bold Q\mid_{D_J^0})$ is given by
the complex
$$
K_J^{\bullet} : \Omega^{\bullet}(J) \to
\oplus_{\# K =1, K\cup J = \emptyset}\Omega^{\bullet}(J \cup K) 
\to
\oplus_{\# K =2, K\cup J = \emptyset}\Omega^{\bullet}(J \cup K) 
\to \cdots
$$
with the filtrations
$$
\align
W_k K_J^{\bullet}:
W_k\Omega^{\bullet}(J)  &\to
\oplus_{\# K =1, K\cup J = \emptyset}W_{k+1}\Omega^{\bullet}(J \cup K)  \\
& \to
\oplus_{\# K =2, K\cup J = \emptyset}W_{k+2}\Omega^{\bullet}(J \cup K) 
\to \cdots, \text{ and } \\
\endalign
$$
$$
\align
F^p K_J^{\bullet}:
F^p\Omega^{\bullet}(J) & \to
\oplus_{\# K =1, K\cup J = \emptyset}F^p\Omega^{\bullet}(J \cup K)  \\
& \to
\oplus_{\# K =2, K\cup J = \emptyset}F^p\Omega^{\bullet}(J \cup K) 
\to \cdots. \\
\endalign
$$
We define sub-complex $K^{\bullet}_{J,i}$, 
$W_kK^{\bullet}_{J,i}$ and $F^pK^{\bullet}_{J,i}$ by
$$
\align
K_{J,i}^{\bullet}:
W_i\Omega^{\bullet}(J) & \to
\oplus_{\# K =1, K\cup J = \emptyset}W_{i}\Omega^{\bullet}(J \cup K)  \\
& \to
\oplus_{\# K =2, K\cup J = \emptyset}W_{i}\Omega^{\bullet}(J \cup K) 
\to \cdots, \\ 
\endalign
$$
$W_kK_{J,i}^{\bullet} = K^{\bullet}_{J,i} \cap W_k K_J^{\bullet}$
and
$F^pK_{J,i}^{\bullet} = K^{\bullet}_{J,i} \cap F^p K_J^{\bullet}$.
\proclaim{Proposition 3.5.1}
\roster
\item
The sequence of morphism 
$K_{J,i}^{\bullet} \to K_{J}^{\bullet}$,
$K_{J, i-1}^{\bullet} \to K_{J,i}^{\bullet}$ are quasi canonical
filtration.
\item
The quotient complex $K_{J,i} / K_{J,i-1}$ is weakly equivalent to
$L_J^{\bullet}\otimes\wedge^i(\bold Q^J)$
$\simeq L_J^{\bullet} \otimes (\oplus_{\# L =i, L \subset J}\bold Q)$,
where
$$
L_J^{\bullet} : \Omega_{D_J}^{\bullet}[-i] \to \oplus_{\# K =1}
\Omega_{D_{J\cup K}}[-i] \to \cdots
$$
\item
Let $\bold Q^J =M_1 \oplus M_2$ be a direct sum decomposition of 
$\bold Q^J$.  Then the direct sum decomposition
$\wedge^i M_1 \oplus M_2 \wedge \wedge^i(\bold Q^J)$ give reise to
a direct sum decomposition 
$L_J^{\bullet} \otimes \wedge^i M_1 \oplus 
L_J^{\bullet}\otimes(M_2 \wedge \wedge^i(\bold Q^J))$
of $K_{J,i}/K_{J,i-1}$.
\endroster
\endproclaim
\demo{Proof}
First we show that $K_{J,i}^{\bullet}$ is a CMHC.  The sheaf $Gr_k^WK_{J,i}$
is given by 
$$
\align
 Gr_k^W K_{J,i}^{\bullet} =& (
Gr_k^W\Omega^{\bullet}(J) \overset{0}\to\to 
\oplus_{\# K =1, K\cup J = \emptyset}Gr_{k+1}^W\Omega^{\bullet}(J \cup K) 
\overset{0}\to\to\cdots \overset{0}\to\to  \\
& \oplus_{\# K =i-k, K\cup J = \emptyset}Gr_{i}^W\Omega^{\bullet}(J \cup K) 
\to 0 ) \\
= &
(\oplus_{\# L =k}
\Omega^{\bullet}_{D_{L\cup J}}[-k] \overset{0}\to\to 
\oplus_{\# K =1,\#L =k+1, K\cup J = \emptyset}
\Omega^{\bullet}_{D_{L\cup J \cup K}}[-k-1] 
\overset{0}\to\to \\
&\cdots \overset{0}\to\to
\oplus_{\# K =i-k,\# L=k+(i-k), K\cup J = \emptyset}
\Omega^{\bullet}_{D_{L\cup J \cup K}}[-i]
\to 0).  \\
\endalign
$$
Therefore, we have 
$$
\align
 \bold H^m(Gr^p_FGr_k^WK_{J,i}) = &
 \oplus_{\# L =k}
H^{m-p}(\Omega^{p-k}_{D_{L\cup J}}) \oplus  \\ 
& \oplus_{\# K =1,\#L =k+1, K\cup J = \emptyset}
H^{p-m-1}(\Omega^{p-k-1}_{D_{L\cup J \cup K}} )\\
& \oplus\cdots \oplus
\oplus_{\# K =i-k,\# L=k+(i-k), K\cup J = \emptyset}
H^{m-i+k-p}(\Omega^{p-i}_{D_{L\cup J \cup K}}).  \\
\endalign
$$
and $Gr_k^W(K_{J,i}^{\bullet})$ is a cohomological Hodge complex.
Since $K_{J,i}^{\bullet}$ for $i << 0$ and $K_{J,i}\simeq K_J$
for $i >> 0$, it is sufficient to prove that 
$\Cal H^p(K_{J, i}^{\bullet}/K_{J,i-1}^{\bullet})=0$
if $p \neq i$.  Since $K^{\bullet}_{J,i}$ is quasi isomorphic to
$$
\tau_i \bold R j_* \bold Q\mid_{D_J} \to
\oplus_{\# K =1, K \cup J = \emptyset}
\tau_i \bold R j_* \bold Q \mid_{D_{J\cup K}} \to \cdots,
$$
the quotient complex $K_{J,i}/K_{J,i-1}$ is quasi isomorphic to
$$
\bold R^i j_* \bold Q\mid_{D_J}[-i] \to
\oplus_{\# K =1, K \cup J = \emptyset}
\bold R^i j_* \bold Q \mid_{D_{J\cup K}}[-i] \to \cdots,
$$ 
and it is quasi isomorphic to
$(j_J)_!(\bold R j_* \bold Q\mid_{D^0_{J}})[-i]$.  This proves the 
proposition.

(2) The associated graded complex of sheaves $K_{J,i}/ K_{J,i-1}$
is isomorphic to 
$$
\oplus_{\# L =i}\Omega^{\bullet}_{L\cup J}[-i]
\to
\oplus_{\# L =i,\# K = 1, K \cap J =\emptyset}
\Omega^{\bullet}_{L\cup J\cup K}[-i]
\to\cdots
$$
and therefore it is a direct sum of $M_L^{\bullet}$, $(\# L =i)$, 
where
$$
M_L^{\bullet} =
(\Omega_{L\cup J}^{\bullet}[-i] \to
\oplus_{\# K =1, K\cap J =1}\Omega^{\bullet}_{L\cup J \cup K}[-i]
\to \cdots).
$$
Moreover, the direct sum is compatible with the filtration $F^{\bullet}$
and $W_{\bullet}$ induced by that of $K_{J,i}^{\bullet}$ and it
is easy to see that these filtrations defines structure of CMHC.
This morphism defines a morphism $\Cone (K_{J,i-1} \to K_{J,i})[1]
\to \oplus_{\# L =i}M^{\bullet}_L$ of CMHC which is a quasi-isomorphism, i.e.
a weak equivalence of CMHC.  (Note that in general if 
$0 \to K \to L \to M \to 0$ is an exact sequence of CMHC, then the
natural map $\Cone (K \to L)[1] \to M$ is a weak equivalence.
Therefore $H^i(X, \Cone (K \to L)) \simeq H^i(X, M)$ is an
isomorphism of mixed Hodge structures. cf. [Dur], \S 2.5 for $X=$ a point.)
We show that $M_L^{\bullet}$ is quasi-isomorphic to $0$ if
$L \nsubset J$.  Let $a \in L -J$ and put $J'=J \cup \{a \}$.
We define complexes $'M_L^{\bullet}$ and $''M_L^{\bullet}$
as
$$
\align
'M_L^{\bullet} & =
(\Omega_{L \cup J'} \to \oplus_{\# K=1, K \cap J' =\emptyset}
\Omega_{L \cup J' \cup K} \to \cdots ) \\
''M_L^{\bullet} & =
(0 \to \oplus_{\# K=1, K \cap J =\emptyset, a \in K}
\Omega_{L \cup J' \cup K} \to 
\oplus_{\# K=2, K \cap J =\emptyset, a \in K}
\Omega_{L \cup J' \cup K} \to \cdots ) \\
\endalign
$$
Then it is easy to see that $M^p_L = 'M^p_L \oplus
''M^p_L$.  The bijection between
$\{ K' \mid \# K' = p , K'\cap J' =\emptyset\}$ and
$\{ K \mid \# K =p+1, K \cap J = \emptyset, a \in K\}$ given by
$K' = K-\{a \}$ gives an isomorphism of complexes
$'M_L^{\bullet} \to ''M_L^{\bullet}$.  As a complex, it is easy
to check that $M^{\bullet}_L = \Cone ('M^{\bullet}_L \to
''M^{\bullet}_L)$ and this proves the exactness of $M_L^{\bullet}$
if $L \nsubset J$.  Moreover, if $L \subset J$, then $M_L^{\bullet}
\simeq L_J^{\bullet}$.  Therefore we get the statement (2) of the
proposition.

The stetement (3) is a direct consequence of (2).
\enddemo

\heading
\S 3.6 Convolution theorem for Hodge theory.
\endheading

First, we compute $\Psi_{V,m,\sigma}(\bold Q)$ using the
result of \S 3.3, \S 3.4 and \S 3.5.  
The normal crossing divisor $V_B$ of $V$ defines a natural stratification
and $V_c$ is expressed as a union of this stratification:
$\cup_{\tau \in K}S_{\tau}^0$ indexed by a set $K$.
The inverse image of $0$ under the natural map $\tilde B(y_i) \to B(y_i)$
is denoted by $V_{c,i}$.
The morphism 
$R^i(\tilde j_M)_*\bold Q \to R^i(\tilde j_{EFB})_*\bold Q$ is 
surjective and identified with a projection for direct sum decomposition
of CMHC's by the expression of 
$R^i(\tilde j_M)_* \bold Q\mid_{\tilde Z_{\sigma}^0}$ in Proposition 3.3.2
and the expression of $\tilde \alpha^*$ in Proposition 3.3.1.(2).  
Therefore $R^i(\tilde j_{EFB})_*\bold Q
\mid_{S_{\tau}^0\times \tilde Z_\sigma^0}$ has a structure of CMHC which
is weakly equivalent to $\wedge^i(\bold Q(-1)^{\tilde r(\sigma, \tau)})$, as
a CMHC, where $\tilde r (\sigma, \tau)
=\dim R^1(\tilde j_{EFB})_*\bold Q \mid_{\bar v}$ for $\bar v \in 
S_{\tau}^0 \times \tilde Z_\sigma^0$ using Proposition 3.5.1.(3) and
Theorem 2.4.9.
We use the same argument again to compute 
$R^i(\tilde j_{EF})_*\bold Q\mid_{S_{\tau}^0\times \tilde Z^0_{\sigma}}$
and we get an isomorphism 
$R^i(j_{EF})_* \bold Q \simeq \wedge^i\bold Q(-1)^{r(\sigma, \tau)}$
as a CMHC,s, where $r(\sigma, \tau) = \sum_{i=k+1}^nl_i +a -n$.
Here we used the same notations $l_{k+1}, \dots , l_n$ and $a$ in \S 3.3.
Therefore, we have an isomorphism
$$
\bold H^*_c(S_{\tau}^0\times\tilde Z_{\sigma}^0,
\bold R^i(\tilde j_{EF})_*\bold Q)
\simeq
\bold H^*_c(S_{\tau}^0, \bold Q)
\otimes\bold H^*_c(\tilde Z_{\sigma}^0, \bold Q)
\otimes \wedge^i(\bold Q(-1)^{r(\sigma, \tau)}),
$$
as a $G_{\sigma, \tau}$-Hodge structure, where
$G_{\sigma, \tau}$ is the stabilizer of 
$S_{\tau} \times \tilde Z^0_{\sigma}$ in $G\times \mu_m$.
As a consequence, we have
$$
[\bold H^*_c(S_{\tau}^0\times\tilde Z_{\sigma}^0,
\bold R^i(\tilde j_{EF})_*\bold Q)]
=
[\bold H^*_c(S_{\tau}^0, \bold Q)]
\otimes[\bold H^*_c(\tilde Z_{\sigma}^0, \bold Q)]
\otimes [\bold Q -\bold Q(-1)]^{r(\sigma, \tau)}.
$$
in $K_{MH}(\bold C, G_{\sigma, \tau})$.
Let $r'(\sigma, \tau)= \sum_{i=k+1}^n(l_i -1)$
and $r''(\sigma, \tau) = a-k$.  Then we have the following equality:
$$
\sum_{[\tau] \in K/G}\Ind_{G_\tau}^G[\bold H^*_c
(S^0_{\tau}, \bold Q)]\otimes 
[\bold Q -\bold Q(-1)]^{r'(\sigma, \tau)}
=
\otimes_{i \in I} [\bold H^*_c(V_{c,i}, \bold Q)] \otimes_{i \notin I}
\phi_{g_i, d_i}(\bold Q),
$$
in $K_{MH}(\bold C, G)$,
where $I = I(\sigma) =\{i \mid \bold e_i \in \sigma \}$.
Therefore, we have
$$
\Psi_{V, m, \sigma}= \otimes_{i \in I}
[\bold H^*_c(V_{c, i}, \bold Q)] \otimes
\otimes_{i \notin I}\phi_{g_i,d_i}(\bold Q) \otimes
[\bold H^*_c(\tilde Z_{\sigma}^0, \bold Q)]
\otimes [\bold Q -\bold Q(-1)]^{r''(\sigma, \tau)}
$$
in $K_{MH}(\bold C, G \times \mu_m)$.
For any $I \subset [1,n]$, we have
$$
\sum_{\{I(\sigma) = I\}/G}
\Ind_{G_\sigma}^G
[\bold H^*_c(\tilde Z^0_{\sigma}, \bold Q)]
\otimes ( \bold Q -\bold Q(-1))^{r''(\sigma, \tau)}
=
\bold H^*_c(\tilde Z^0_I, \bold R \tilde j_*\bold Q),
$$
where $\tilde j : \tilde B(\xi ) -(\tilde D_E \cup \tilde D_F) \to
\tilde B(\xi)$ and $\tilde Z^0_I = \cup_{I(\sigma)=I}\tilde Z_{\sigma}^0$.
As a consequence, we have
$$
\align
& [\bold H^*(\tilde B_\bold d(y)-(f(\xi)\circ g)^{-1}(0), 
k_{\emptyset !}\bold Q)] \\
=&
\sum_{K \times \{\sigma\mid  
I(\sigma) = \emptyset\}/G \times \mu_m}
\Ind_{G_{\sigma, \tau}}^{G \times \mu_m}
[\bold H^*_c(S_{\tau}^0, \bold Q)]
\otimes[\bold H^*_c(\tilde Z_{\sigma}^0, \bold Q)]
\otimes [\bold Q -\bold Q(-1)]^{r(\sigma, \tau)} \\
=&
\otimes_{i=1}^n \phi_{g_i, d_i}(\bold Q) \otimes
[\bold H^*_c(\tilde Z_\emptyset^0, \bold R\tilde j_* \bold Q)], \\
\endalign
$$
where $k_{\emptyset}$ and $\tilde B^0_d(y)$ is defined at the end of \S 3.2.
In the same way, for a subset $I$ of $[1,n]$, we have
$$
[\bold H^*(\tilde B_{I,\bold d}(y)-(f_I(\xi)\circ g)^{-1}(0), 
k_{I !}\bold Q)]
=
\otimes_{i\notin I} \phi_{g_i, d_i}(\bold Q) \otimes
[\bold H^*_c(\tilde Z_I^0, \bold R\tilde j_* \bold Q)].
$$
Using the equality (3.2.2), we have 
$$
\align
&[\bold H^*(\tilde B_\bold d(y)-(f(\xi)\circ g)^{-1}(0), \bold Q)] \\
=&
\sum_{I \subset [1,n]}
[\bold H^*(\tilde B_{I,\bold d}(y)-(f_I(\xi)\circ g)^{-1}(0), 
k_{I !}\bold Q)] \\
=&\sum_{I \subset [1, n]}
\otimes_{i\notin I} \phi_{g_i, d_i}(\bold Q) \otimes
[\bold H^*_c(\tilde Z_I^0, \bold R\tilde j_{I*} \bold Q)] \\
=& \sum_{I \subset [1, n]}(-1)^{\# I}
\otimes_{i\notin I} \phi_{g_i, d_i}(\bold Q) \otimes
\otimes_{i\in I} (\phi_{g_i, d_i}(\bold Q)-[\bold Q]) \otimes
[\bold H^*_c(\tilde Z_I, \bold R\tilde j_{I*} \bold Q)] \\
= & 
\sum_{I \subset [1, n]}(-1)^{\# I}
\otimes_{i\notin I} \phi_{g_i, d_i}(\bold Q) \otimes
\otimes_{i\in I} (\phi_{g_i, d_i}(\bold Q)-[\bold Q]) \\
&\qquad \qquad\otimes
([\bold H^*_c(\tilde Z_I, \bold R\tilde j_{I*} \bold Q)]-[\bold Q-
\bold Q(-1)] \\
& + [\bold Q - \bold Q (-1)] \\
\endalign
$$
we have
$$
\align
\tilde\Phi_{f\circ g, m}(\bold Q)
& =(\frac{1}{[\bold Q -\bold Q(-1)]}
[\bold H_c^*(\tilde B_\bold d(y)-(f(\xi)\circ g)^{-1}(0), \bold Q)] )
-\bold Q \\
& =
\sum_{I \subset [1, n]}(-1)^{\# I}
\otimes_{i\notin I} \phi_{g_i, d_i}(\bold Q) \otimes
\otimes_{i\in I} \Phi_{g_i, d_i}(\bold Q) \otimes
\tilde\Phi_{f_I,m}(\bold Q). \\
\endalign
$$

By taking $G$-invariant part of the above equality, 
we have the following theorem.
Note that the action of $G$ on 
$\bold H_c^*(\tilde Z_I^0, \bold R \tilde j_* \bold Q)$ factors
through $G_I$.
\proclaim{Theorem 3.6.1 (Convolution Theorem for Hodge structures)}
$$
\Phi_{f\circ g, m}(\bold Q) =
\sum_{I \subset [1, n]}(-1)^{\# I}
(\otimes_{i \notin I}
\phi_{g_i, d_i}(\bold Q) 
\otimes\tilde\Phi_{f_I,m}(\bold Q))^{G_I}
\otimes
\otimes_{i \in I}\Phi_{g_i, d_i}(\bold Q)^{\mu_{d_i}}
$$
\endproclaim
By specializing to the case where $n=2$ and $f=x_1 + x_2$, we
have the following corollary.
\proclaim{Corollary 3.6.2 (Thom-Sebastiani theorem for Hodge structures)}
Let $g_1$ and $g_2$ be germs of holomorphic functions on
$\bold C^{n_1}$ and $\bold C^{n_2}$ at $0$ with isolated singularities.
Let $d_1$ and $d_2$ be the
exponents of $g_1$ and $g_2$ respectively and $m = \lcm (d_1, d_2)$.
Then we have
$$
\align
\Phi_{g_1+g_2,m} = &
(\Phi_{g_1,d_1}\otimes \Phi_{g_2,d_2} \otimes
\tilde\Phi_{\xi_1^{d_1}+\xi_2^{d_2}})^{\mu_{d_1} \times \mu_{d_2}} \\
& -\Phi_{g_1,m,\neq 1}\otimes\Phi_{g_2,m,\neq 1} 
 +\Phi_{g_1,m}\otimes\Phi_{g_2,m},\\
\endalign
$$
where $\Phi_{g_1+g_2, m}=\Phi_{g_1+g_2, m}(\bold Q)$ e.t.c.
\endproclaim
\demo{Proof}
By the Theorem 3.6.1, we have
$$
\align
\Phi_{g_1+g_2,m} = &
(\phi_{g_1,d_1}\otimes \phi_{g_2,d_2} \otimes
\tilde\Phi_{\xi_1^{d_1}+\xi_2^{d_2},m})^{\mu_{d_1} \times \mu_{d_2}} 
\tag{3.6.1}\\
& -(\phi_{g_1,d_1}\otimes 
\tilde\Phi_{\xi_1^{d_1},m})^{\mu_{d_1}}
\otimes \Phi_{g_2,d_2}^{\mu_{d_2}} \\
& -(\phi_{g_2,d_2}\otimes 
\tilde\Phi_{\xi_2^{d_2},m})^{\mu_{d_2}}
\otimes \Phi_{g_1,d_1}^{\mu_{d_1}} \\
& +\Phi_{g_1,d_1}^{\mu_{d_1}}\otimes\Phi_{g_2,d_2}^{\mu_{d_2}}. \\
\endalign
$$
On the other hand, we have
$$
\tilde\Phi_{\xi_1^{d_1},m} = [\bold Q[\mu_{d_1}]]-[\bold Q]
\tag{3.6.2}
$$
in $K_{MH}(\bold C, \mu_{d_1}\times \mu_m)$,
where $\mu_m$ acts on $\mu_{d_1}$ through its quotient.
Therefore we have
$$
( \phi_{g_1,d_1} \otimes \tilde\Phi_{\xi^{d_1},m})^{\mu_{d_1}}
=\phi_{g_1,d_1, \neq 1}=
\Phi_{g_1,d_1, \neq 1}.
$$
If
$
(\tilde\Phi_{\xi_1^{d_1} + \xi_2^{d_2},m}\otimes \bold C)(\chi_1, \chi_2)
\neq 0
$
for a charcters $\chi_1$ and $\chi_2$ of $\mu_{d_1}$ and $\mu_{d_2}$,
then $\chi_1 \neq 1$ and $\chi_2 \neq 1$.
Therefore we have
$$
(\phi_{g_1,d_1}\otimes \phi_{g_2,d_2} \otimes
\tilde\Phi_{\xi_1^{d_1}+\xi_2^{d_2},m})^{\mu_{d_1} \times \mu_{d_2}} 
=
(\Phi_{g_1,d_1}\otimes \Phi_{g_2,d_2} \otimes
\tilde\Phi_{\xi_1^{d_1}+\xi_2^{d_2},m})^{\mu_{d_1} \times \mu_{d_2}} 
\tag{3.6.3}
$$
By (3.6.1), (3.6.2) and (3.6.3), we have the corollary.
\enddemo

\enddocument